\pdfoutput=1
\documentclass[12pt]{amsart}

\usepackage[style=numeric,backend=biber,sorting=nyt]{biblatex}
\usepackage{adjustbox, afterpage, amsbsy, amsmath, amssymb, amsthm,
  caption, diagbox, enumitem, fixmath, graphicx, mathtools,
  mdframed, pigpen, pgfkeys, subcaption,
  xspace,url}   

\usepackage[dvipsnames]{xcolor}

\usepackage[color,matrix,arrow,2cell,all, rotate,pdf]{xy}
\UseAllTwocells
\CompileMatrices

\newdir{ >}{!/6pt/@{ }*:(1,0.2)@_{>}*:(1,-0.2)@^{>}}

\newcommand{\BiTop}{\mbox{$\pmb{BiTop}$}\xspace}

\newcommand{\Cat}{\mbox{$\pmb{Cat}$}\xspace}

\newcommand{\CRing}{\mbox{$\pmb{CRing}$}\xspace}

\newcommand{\Epi}[1]{\mbox{$\mathrm{Epi}(#1)$}\xspace}

\newcommand{\ExtEpi}[1]{\mbox{$\mathrm{ExtEpi}(#1)$}\xspace}
\newcommand{\ExtMon}[1]{\mbox{$\mathrm{ExtMon}(#1)$}\xspace}

\newcommand{\Fil}[1]{\mbox{$\pmb{Fil}(#1)$}\xspace}

\newcommand{\Int}[2]{\mbox{$#1[#2]$}\xspace}

\newcommand{\InvRAdjt}[1]{\mbox{${#1}_{\pmb{ppj}}$}\xspace}

\newcommand{\Iso}[1]{\mbox{$\mathrm{Iso}(#1)$}\xspace}

\newcommand{\Loc}{\mbox{$\pmb{Loc}$}\xspace}

\newcommand{\Meas}{\mbox{$\pmb{Meas}$}\xspace}

\newcommand{\Mono}[1]{\mbox{$\mathrm{Mon}(#1)$}\xspace}

\newcommand{\NHD}[1]{\mbox{$\pmb{Nbd}[#1]$}\xspace} % the category
\newcommand{\nhd}[1]{\mbox{$\pmb{nbd}[#1]$}\xspace} % neighbourhood structures

\newcommand{\pfs}[1]{\mbox{$\pmb{Pfs}[#1]$}\xspace}

\newcommand{\pNHD}[1]{\mbox{$\pmb{pNbd}[#1]$}\xspace} % pre-nbds on an
\newcommand{\pnhd}[1]{\mbox{$\pmb{pnbd}[#1]$}\xspace} % category

\newcommand{\PPJ}[1]{\mbox{$#1_{\pmb{ppj}}$}\xspace}

\newcommand{\pre}[1]{\mbox{$\pmb{pre}[#1]$}\xspace}

\newcommand{\RegEpi}[1]{\mbox{$\mathrm{RegEpi}(#1)$}\xspace}

\newcommand{\Set}{\mbox{$\pmb{Set}$}\xspace}

\newcommand{\slice}[2]{\mbox{$( #1 \downarrow\, #2)$}\xspace}

\newcommand{\Sub}[2]{\mbox{$\pmb{Sub}_{#2}(#1)$}\xspace}

\newcommand{\tnhd}[1]{\mbox{$\pmb{top}[#1]$}\xspace}
\newcommand{\Top}{\mbox{$\pmb{Top}$}\xspace}

\newcommand{\TychObj}[1]{\mbox{$A$-$\pmb{Tych}$}\xspace}

\newcommand{\wNHD}[1]{\mbox{$\pmb{wNbd}[#1]$}\xspace} % category of
\newcommand{\wnhd}[1]{\mbox{$\pmb{wnbd}[#1]$}\xspace} % weak
\newcommand{\adjt}[2]{\mbox{$ {#1} \dashv {#2}$}}

\newcommand{\adjtpair}[5]{\mbox{$\xymatrix@!C=#5{ {#3} 
\ar@/^2ex/[r]^{#1}  & {#4} \ar@/^2ex/[l]^{#2} }$}}
\newcommand{\ADJTSIT}[7]{\mbox{$\Bigl\langle\xymatrix{ {#3}
\ar@<0.9ex>[r]^{#1} \ar@{}[r]|{\bot} \ar@<-0.9ex>@{<-}[r]_{#2} & {#4} };
\mathrm{unit:}\xspace \id{#3} \xrightarrow{#5} {#2}\circ {#1},
\mathrm{counit:}\xspace {#1}\circ {#2} \xrightarrow{#6}
\id{#4}; \mathrm{adjt. iso:}\xspace {#7}\Bigr\rangle$}} % first adjt
\newcommand{\Adjtsit}[6]{\mbox{$\Bigl\langle\xymatrix{ {#3}
\ar@<0.9ex>[r]^{#1} \ar@{}[r]|{\bot} \ar@<-0.9ex>@{<-}[r]_{#2} & {#4} };
\mathrm{unit:}\xspace \id{#3} \xrightarrow{#5} {#2}\circ {#1},
\mathrm{counit:}\xspace {#1}\circ {#2} \xrightarrow{#6} 
\id{#4}\Bigr\rangle$}} % first adjt diagram, follows by unit, co-unit

\newcommand{\adjtsit}[4]{\mbox{$\xymatrix@!{ {#3} \ar@<0.9ex>[r]^{#1} \ar@{}[r]|{\bot} \ar@<-0.9ex>@{<-}[r]_{#2} & {#4} }$}}

\newcommand{\adjtequiv}[4]{\mbox{$\xymatrix@!{ {#3} \ar@<0.9ex>[r]^{#1} \ar@{}[r]|{\top\text{\hspace{-5pt}}\bot} \ar@<-0.9ex>@{<-}[r]_{#2} & {#4} }$}}

\newcommand{\ArGrFib}[7]{\mbox{$\xymatrix{ {#4} \ar@{->}[dd]_{#2} & {#6} \ar@{->}[dd]^{#3} \\ {} \ar@{->}[r]^{#1} & \\ {#5} & {#7} }$}}
\newcommand{\Arr}[3]{\mbox{$#2 \xrightarrow{#1} #3$}}

\newcommand{\arrp}[4]{\mbox{$\xymatrix@!{ {#3} \ar@<3pt>[r]^{#1}
      \ar@<-3pt>[r]_{#2} & {#4} }$}}
\newcommand{\arrpGrFib}[8]{\mbox{$\xymatrix{ {#5} \ar@{->}[dd]_{#3} &
{#7} \ar@{->}[dd]^{#4} \\ {} \ar@<1ex>[r]^{#1} \ar@<-1ex>[r]_{#2} & \\
{#6} & {#8} }$}}

\newcommand{\atleast}[1]{\mbox{$\uparrow #1$}}
\newcommand{\Bb}[1]{\mbox{${\mathbb #1}$}}

\newcommand{\CAL}[1]{\mbox{${\mathcal #1}$}}

\newcommand{\comp}[2]{\mbox{${#1}{\circ}{#2}$}}

\newcommand{\commsq}[8]{\mbox{$\xymatrix@!{ {#4} \ar@{->}[r]^{#8} \ar@{->}[d]_{#7} \ar@{}[dr]|{\#} & {#1} \ar@{->}[d]^{#5} \\ {#3} \ar@{->}[r]_{#6} & {#2} }$}}

\newcommand{\dn}[1]{\mbox{${#1}^{\downarrow}$}}

\newcommand{\down}[2]{\mbox{$#1 \downarrow #2$}}

\newcommand{\epi}[3]{\mbox{$\xymatrix@!{ {#2} \ar@{->>}[r]^{#1} & {#3} }$}}

\newcommand{\eulf}[1]{\mbox{${\mathfrak #1}$}}
\newcommand{\fact}[2]{$(#1, #2)$-factorisation}

\newcommand{\finv}[2]{\mbox{${#1}^{^{\mathbf{-1}}}{#2}$}}
\newcommand{\ffinv}[3]{\mbox{$\bigl({#1}^{-1}{#2}\bigr)^{-1}{#3}$}\xspace}
\newcommand{\fInv}[2]{\mbox{$\forall_{_{\hspace{-3pt}#1}}{#2}$}\xspace}

\newcommand{\galstr}[8]{\mbox{$\Bigl\langle\Bigl\langle\xymatrix{ {#3}  
\ar@/^1.56ex/[r]^{#1} \ar@{}[r]|{\bot} & {#4} \ar@/^1.56ex/[l]^{#2} }; 
#5, #6; #7, #8\Bigr\rangle\Bigr\rangle$}}

\newcommand{\ID}[1]{\mbox{$\xymatrix{ {#1} \ar@2{-}[r]^{} & {#1} }$}}
\newcommand{\id}[1]{\mbox{$\mathbf{1}_{#1}$}}

\newcommand{\img}[2]{\mbox{$\exists_{_{#1}}{#2}$}}
\newcommand{\imgfil}[2]{\mbox{$\overset{\rightarrow}{#1}{#2}$}\xspace}
\newcommand{\Imgfil}[2]{\mbox{$\overset{\sqcap}{#1}{#2}$}\xspace}

\newcommand{\intr}[2]{\mbox{{\rm int}$_{#1}#2$}\xspace}

\newcommand{\invfil}[2]{\mbox{$\overset{\leftarrow}{#1}{#2}$}\xspace}
\newcommand{\Invfil}[2]{\mbox{$\overset{\sqcup}{#1}{#2}$}\xspace}
\newcommand{\InvRel}[3]{\mbox{$#1^{-1}\bigl(\xymatrix{ {#2} 
\ar@/^2ex/[r] \ar@/_2ex/[r] & {#3} }\bigr)$}\xspace}
\newcommand{\ISO}[3]{\mbox{$\xymatrix{ {#2} \ar@{->}[r]^{#1}
      \ar@{}[r]_{\text{iso}} & {#3} }$}}
\newcommand{\iso}[3]{\mbox{$\xymatrix@!{ {#2} \ar@{->}[r]^{#1}
\ar@{}[r]_{\text{\tiny iso}} & {#3}}$}}

\newcommand{\Kerp}[1]{\mbox{$\mathrm{Kerp}[#1]$}\xspace}

\newcommand{\limit}[2]{\mbox{$\underset{\longleftarrow}{\lim}_{{#2}} #1$}}

\newcommand{\low}[2]{\mbox{${#1}_{{#2}}$}}

\newcommand{\mon}[3]{\mbox{$\xymatrix{ {#2} \ar@{ >->}[r]^{#1} & {#3}}$}}

\newcommand{\Nattrn}[5]{\mbox{$\xymatrix@!{ {#4} \rtwocell^{#2}_{#3}{\,#1} & {#5}}$}\xspace}
\newcommand{\NattrnGrFib}[9]{$\xymatrix@!C=12ex{ {#6} \ar@{->}[dd]_{#4} & {#8} \ar@{->}[dd]^{#5} \\ {} \rtwocell^{#3}_{#2}{!\,\,\,#1} & {} \\ {#7} & {@9} }$}
\newcommand{\nattrn}[3]{\mbox{$\xymatrix{ {#2} \ar@2{->}[r]^{#1} & {#3} }$}\xspace}

\newcommand{\normalto}[3]{\mbox{$#1 \vartriangleleft \bigl(\xymatrix{#2
\ar@/^1ex/[r] \ar@/_1ex/[r] & {#3} }\bigr)$}\xspace}

\newcommand{\opair}[2]{\mbox{$(#1, #2)$}}
\newcommand{\opp}[1]{\mbox{$#1^{\mathrm{op}}$}}

\newcommand{\overfib}[2]{\mbox{$\xymatrix{ {#1} \ar@{|->}[r] & {#2} }$}\xspace}
\newcommand{\Overfib}[3]{\mbox{$\xymatrix{ {#2} \ar@{|->}[r]^{#1} & {#3}
}$}\xspace}
\newcommand{\overfibp}[6]{\mbox{$\xymatrix{ {#3} \ar[r]^{#1} \ar@{|->}[d] & 
{#4} \ar@{|->}[d] \\ {#5} \ar[r]_{#2} & {#6} }$}\xspace}
\newcommand{\Overfibp}[7]{\mbox{$\xymatrix{ {#4} \ar[r]^{#2} 
\ar@{|->}[d]_{#1} & {#5} \ar@{|->}[d]^{#1} \\ {#6} \ar[r]_{#3} & {#7}
}$}\xspace}

\newcommand{\pow}[2]{\mbox{${#1}^{{#2}}$}}

\newcommand{\Pulb}[7]{\mbox{$\xymatrix@!{ {#1 \times_{#3} #2}
\ar@{->}[r]^{#4} \ar@{->}[d]_{#5} \ar@{}[dr]|(0.096){\text{\LARGE{$\lrcorner$}}}& {#2} \ar@{->}[d]^{#6} \\ {#1} \ar@{->}[r]_{#7} & {#3} }$}}
\newcommand{\pulb}[8]{\mbox{$\xymatrix@!{ {#1} \ar@{->}[r]^{#5}
\ar@{->}[d]_{#6} \ar@{}[dr]|(0.096){\text{\LARGE{$\lrcorner$}}} & {#2}
\ar@{->}[d]^{#7} \\ {#3} \ar@{->}[r]_{#8} & {#4}  }$}}
\newcommand{\pulbint}[7]{\mbox{$\xymatrix@!{ {#1 \cap #3} \ar@{ >->}[r]^{#7} \ar@{ >->}[d]_{#6} \ar@{ >->}[dr]|{#4 \cap #5} & {#1} \ar@{ >->}[d]^{#4} \\ {#3} \ar@{ >->}[r]_{#5} & {#2} }$}}
\newcommand{\Pulbmon}[8]{\mbox{$\xymatrix@!{ {#4 = {#7}^{\leftarrow}(#1)} \ar@{->}[r]^{#8} \ar@{ >->}[d]_{#7} & {#1} \ar@{ >->}[d]^{#5} \\ {#3} \ar@{->}[r]_{#6} & {#2} }$}}
\newcommand{\pulbmon}[8]{\mbox{$\xymatrix@!{ {#4} \ar@{->}[r]^{#8} \ar@{
>->}[d]_{#7} \ar@{}[dr]|{\text{pullback}} & {#1} \ar@{ >->}[d]^{#5} \\
{#3} \ar@{->}[r]_{#6} & {#2} }$}} 

\newcommand{\Push}[7]{\mbox{$\xymatrix@!{ {#3} \ar@{->}[r]^{#4}
\ar@{->}[d]_{#5} \ar@{}[dr]|(0.906){\text{\large{$\ulcorner$}}} & {#1}
\ar@{->}[d]^{#6} \\ {#2} \ar@{->}[r]_{#7} & {{#1} +_{#3} {#2}} }$}}
\newcommand{\push}[8]{\mbox{$\xymatrix@!{ {#1} \ar@{->}[r]^{#5}
\ar@{->}[d]_{#6} \ar@{}[dr]|(0.906){\text{\large{$\ulcorner$}}} & {#2}
\ar@{->}[d]^{#7} \\ {#3} \ar@{->}[r]_{#8} & {#4} }$}}

\newcommand{\rel}[5]{$\xymatrix{ & {#3} \ar@{->}[dl]_{#1}
\ar@{->}[dr]^{#2} & \\ {#4} & & {#5} }$}
\newcommand{\relon}[6]{\mbox{$\xymatrix@!C=#6ex{ {#4} \ar@/^2ex/[r]^{#1}    
\ar@/_2ex/[r]_{#2} \ar@{<-}[r]|{#3} & {#5} }$}\xspace}  
\newcommand{\Relon}[2]{\mbox{$\xymatrix{ {#1} \ar@/^1ex/[r]   
\ar@/_1ex/[r] & {#2} }$}\xspace} 
\newcommand{\relation}[3]{\mbox{$\xymatrix{ {#2} \ar@{~}[r]^{#1} & {#3} }$}}

\newcommand{\rest}[2]{\mbox{${#1}\bigr|_{#2}$}}

\newcommand{\seq}[3]{\mbox{$\bigl(#1_{#2}\bigr)_{#2 \in #3}$}}

\newcommand{\ses}[5]{\mbox{$\xymatrix{ {\zeroo} \ar@{->}[r]^{} & {#3} 
\ar@{->}[r]^{#1} & {#4} \ar@{->}[r]^{#2} & {#5} \ar@{->}[r]^{} &
{\zeroo} }$}}

\newcommand{\slopen}[1]{\mbox{$\mathfrak{o}({#1})$}\xspace}

\newcommand{\splepi}[4]{\mbox{$\xymatrix{ {#3} \ar@<0.6ex>[r]^{#1} 
\ar@{<-}@<-0.6ex>[r]_{#2} & {#4} }$}\xspace}
\newcommand{\Splepi}[3]{\mbox{$\xymatrix{ {#2} \ar@<0.6ex>[r]^{#1} 
\ar@{<-}@<-0.6ex>[r]_{s_{#1}} & {#3} }$}\xspace}
\newcommand{\splrel}[7]{$\xymatrix{ & {#5} \ar@{->}[dl]_{#1}
\ar@{->}[dr]^{#2} & \\ {#6} \ar@/^8ex/[ur]^{#3} & & {#7}
\ar@/_8ex/[ul]_{#4} }$}

\newcommand{\tfae}{The following are equivalent\xspace}

\newcommand{\up}[1]{\mbox{${#1}^{\uparrow}$}}

\newcommand{\zeroo}{\mbox{$\mathbf{0}$}}

\newenvironment{Cor}{\begin{cor}}{\end{cor}}

\newenvironment{Df}{\begin{df}}{\end{df}}

\newenvironment{Thm}{\begin{thm}}{\end{thm}}
\allowdisplaybreaks
\usepackage[left=2.8cm,top=2.4cm,right=2.8cm,bottom=2.4cm,foot=0.96cm,head=0.96cm]{geometry}  

\swapnumbers
\theoremstyle{plain}

\newtheorem*{thm}{Theorem}
\newtheorem*{lem}{Lemma}
\newtheorem*{prop}{Proposition}
\newtheorem*{cor}{Corollary}

\theoremstyle{definition}

\newtheorem*{df}{Definition}
\newtheorem*{ex}{Example}

\newcommand{\cannbd}[1]{\mbox{${#1}_{\mathbold{w}}$}}
\newcommand{\cansnbd}[1]{\mbox{${#1}_{\mathbold{n}}$}}
\newcommand{\cantop}[1]{\mbox{${#1}_{\mathbold{t}}$}}
\newcommand{\itop}[1]{\mbox{$\mathbold{top}[#1]$}\xspace}
\newcommand{\mo}{\mbox{$\mathfrak{O}$}}

\title{Internal Neighbourhood Structures}
\author{Partha Pratim Ghosh}
\address{Department of Mathematical Sciences \\
Unisa Science Campus \\
corner of Christiaan de Wet \& Pioneer Avenue \\
Florida 1709 \\
Johannesburg, Gauteng \\
South Africa}
\email{ghoshpp@unisa.ac.za}

\begin{document}

\begin{abstract}
 The main aim of this paper is to provide a description of neighbourhood
 operators in finitely complete categories with finite coproducts and a
 proper factorisation system such that the semilattice of {\em
 admissible} subobjects make a distributive complete lattice. The
 equivalence between neighbourhoods, Kuratowski interior operators and
 {\em pseudo-frame sets} is proved. Furthermore the categories of
 internal neighbourhoods is shown to be topological. Regular
 epimorphisms of categories of neighbourhoods are described and
 conditions ensuring hereditary regular epimorphisms are
 probed. It is shown the category of internal neighbourhoods of
 topological spaces is the category of bitopological spaces, while in
 the category of locales every locale comes equipped with a
 natural internal topology. 
 \end{abstract}


\maketitle
%\tableofcontents

\section{Introduction}
\label{Introduction}

The introduction in \cite{DikranjanGiuli1987} and
\cite{DikranjanGiuliTholen1989} of {\em categorical closure operators}
led to systematic study of topological properties in general
categories. The theory of categorical closure operators was subsequently
developed by many authors, for instance in \cite{ClementinoGiuliTholen1996},
\cite{CastelliniGiuli2001}, \cite{MMCEGWT2004}, \cite{Slapal2005} and
the references therein. A concise treatment of this development is
available in the self-contained monograph \cite{DikranjanTholen1995}, as
well as in the later published book \cite{Castellini2003}.   

Closure operators give rise to the notion of {\em neighbourhood
operators} on a category. Categorical neighbourhoods have been treated
in \cite{Slapal2001}, \cite{Slapal2008}, \cite{GiuliSlapal2009},
\cite{HolgateSlapal2011}, \cite{Razafindrakato2012}, \cite{Slapal2012},
\cite{HolgateIragiRazafindrakatos2016},\cite{HolgateRazafindrakato2017}. Since
neighbourhoods are required for the study of convergence, investigation
of neighbourhood structures is important in its own right apart from
being a consequence of the notion of a closure operator.

The purpose of this paper is to show that the notion of a neighbourhood
on an object of a category can be provided with minimal
assumptions. In this paper we assume \Bb{A} to be a finitely complete
category with finite coproducts equipped with a proper
\fact{\mathsf{E}}{\mathsf{M}}  system, such that for each object $X$,
the set \Sub{X}{\mathsf{M}} of $\mathsf{M}$-subobjects of $X$ is a
distributive complete lattice (see page \pageref{basicassumption}). 

Obviously, one can extract more information when stronger
properties of the lattice of admissible subobjects is assumed --- for
instance, Boolean algebra. However, 
the usual correspondence between neighbourhoods, Kuratowski interiors,
and {\em pseudo-frame sets} (see Definition \ref{pseudoframesets-df}, page
\pageref{pseudoframesets-df}) can be obtained in
this general setup (see Theorem
\ref{Kuratowskiisretractofpnbdwithintopen}, page
\pageref{Kuratowskiisretractofpnbdwithintopen} and Theorem
\ref{pfs=nhd}, page \pageref{pfs=nhd} for details). In the remainder of
this introductory section we shall highlight the major results that have
been obtained in this general context.

The notion of a neighbourhood has three layers. The description of these layers involve a filter {\em on}
an object. A filter {\em on} an object $X$ is a filter {\em in} the lattice
\Sub{X}{\mathsf{M}} of admissible subobjects of $X$. As soon as the
lattice of admissible subobjects is not atomic, a collection of
neighbourhoods of an admissible subobject --- which is a
filter, becomes an order reversing assignment from the lattice of
admissible subobjects of the object to the ordered set of filters on the
object. This is the first layer in the definition of a neighbourhood,
herein called {\em preneighbourhoods} (see Definition
\ref{neighbourhoodstructures-df}\ref{preneighbourhood-df}, page
\pageref{preneighbourhood-df}). The second layer in the definition of a
neighbourhood is its {\em interpolability} --- given a neighbourhood $N$
of a subobject $P$, it must be possible to obtain a neighbourhood $N'$
of $P$ such that $N$ is a neighbourhood of $N'$. These neighbourhoods
are called {\em weak neighbourhoods} (see Definition
\ref{neighbourhoodstructures-df}\ref{weakneighbourhood-df}, page
\pageref{weakneighbourhood-df}). Finally come the {\em neighbourhoods}
(see Definition \ref{neighbourhoodstructures-df}\ref{neighbourhood-df},
page \pageref{neighbourhood-df}), a collection of which preserve
arbitrary meets.  

In the special case when the lattice \Sub{X}{\mathsf{M}} of admissible
subobjects of the object $X$ is a frame, neighbourhoods (in the sense of
Definition \ref{neighbourhoodstructures-df}\ref{neighbourhood-df},
hereafter) become topologies on $X$. Furthermore, when the lattice of
admissible subobjects is atomic, the topologies are provided by
prescribing the neighbourhoods of the atoms. In general, a {\em
topology} on $X$ (see Definition \ref{internaltopology-df}, page
\pageref{internaltopology-df}) is a special collection of
neighbourhoods, the set of open sets (see equation \eqref{opensets},
page \pageref{opensets}) is a frame in the order induced from
the lattice \Sub{X}{\mathsf{M}}.    

An object along with a preneighbourhood, or a weak neighbourhood, or a
neighbourhood, or a topology is said to be an {\em internal
preneighbourhood space}, or an {\em internal weak neighbourhood space},
or an {\em internal neighbourhood space}, or an {\em internal
topological space}, respectively. To define the notion of {\em
continuous maps} for these {\em spaces}, one requires the notion of a
preimage. This is achieved from the proper factorisation available on
\Bb{A}. A {\em preneighbourhood morphism} is a morphism $f$ of \Bb{A}
with the property: if $U$ be a preneighbourhood of an admissible
subobject $P$ of the codomain of $f$ then \finv{f}{U} is a
preneighbourhood of \finv{f}{P} (see Definition
\ref{prenbdmorphisms-df}, page \pageref{prenbdmorphisms-df}).    

The collection of internal preneighbourhood spaces of
a category \Bb{A} along with the preneighbourhood morphisms make the category
\pNHD{\Bb{A}}. The full subcategory of internal weak neighbourhood
spaces make the category \wNHD{\Bb{A}}. \wNHD{\Bb{A}} is bireflective in
\pNHD{\Bb{A}} (see 
\ref{wnbdbireflectiveinpretop}, page \pageref{wnbdbireflectiveinpretop} 
and Theorem \ref{nbdreflinpnbd}, page
\pageref{nbdreflinpnbd}). Furthermore, both \pNHD{\Bb{A}} and
\wNHD{\Bb{A}} are topological over \Bb{A} (see Theorem
\ref{topologicityresults}(\ref{pnbdtopoverbase} \&
\ref{wnbdtopoverbase}), page \pageref{topologicityresults}). The
category \pNHD{\Set} of internal preneighbourhood spaces of \Set is
equivalent to the category $\pmb{preTop}$ of pretopological spaces. The
category $\pmb{preTop}$ of pretopological spaces is investigated in
\cite{Kent1969}, \cite{BentleyHerrlichLowen1991} and
\cite{HerrlichLowenSchwarz1991}.  

The morphisms are restricted in the category of internal neighbourhood
spaces. This is suggested from existence of the largest
neighbourhood smaller than a weak neighbourhood in Theorem
\ref{snbdreflinwnbd} (see page \pageref{snbdreflinwnbd}).  The internal
neighbourhood spaces along with preneighbourhood morphisms $f$ for which
the preimage \finv{f}{} preserve arbitrary joins constitute the
subcategory \NHD{\Bb{A}} of internal neighbourhood spaces. Since
morphisms of neighbourhoods are restricted to those whose preimage
preserve joins, it is topological over \PPJ{\Bb{A}} (see Theorem
\ref{topologicityresults}\ref{nbdtopoverppj}, page
\pageref{nbdtopoverppj}), where \PPJ{\Bb{A}} is the non-full subcategory
of \Bb{A} having same objects as \Bb{A} and precisely those morphisms of
\Bb{A} whose preimages preserve joins. Furthermore, \NHD{\Bb{A}} is
bireflective in \PPJ{\wNHD{\Bb{A}}} (see \ref{nbdbireflectiveinwnbd-ppj}, page
\pageref{nbdbireflectiveinwnbd-ppj} and Theorem \ref{snbdreflinwnbd},
page \pageref{snbdreflinwnbd}).   

The full subcategory of \NHD{\Bb{A}} consisting of internal topological
spaces is \Int{\Top}{\Bb{A}}. \Int{\Top}{\Bb{A}} is reflective in
\NHD{\Bb{A}} if and only if \Int{\Top}{\Bb{A}} is topological over
\PPJ{\Bb{A}}  if and only if  each object has a largest internal
topology (see Theorem \ref{largesttop-alt}, page
\pageref{largesttop-alt}). In categories, as in \Set, where each lattice
of admissible subobjects is a 
frame, every object has a largest internal topology.
However, internal topologies are interesting in many situations beyond
\Set. For instance, \Loc is a category in which not every lattice of
admissible subobjects is a frame and the preimage of every
morphism preserve only finite joins. Yet, as shown in the papers
\cite{DubeIghedo2016} \& \cite{DubeIghedo2016a}, the open sublocales do
provide a natural way to define an internal topology
\Arr{\mathit{o}_X}{\opp{\Sub{X}{RegMon}}}{\Fil{X}} (see equation
\eqref{naturaltopologyonlocales-eq}, page
\pageref{naturaltopologyonlocales-eq}) on each locale $X$. In fact, as 
observed in Theorem \ref{Thembanbdrtinv} (see page
\pageref{Thembanbdrtinv}), the assignment $X \mapsto (X, \mathit{o}_X)$
on a locale defines a right inverse to the forgetful functor
\Arr{U}{\pNHD{\Loc}}{\Loc}. 

The lattice of admissible subobjects play an important role in the
development of this paper. Further, apart from the category \Set of sets
and functions there are several categories which satisfy the
basic assumption of this paper (see page \pageref{basicassumption}). The
following is a list of such instances, apart from \Set:
\label{contextlist}
\begin{enumerate}[itemsep=1.2ex,label=(\roman*),align=left]
 \item \Top satisfy the conditions (see \S \ref{topologies}, page
       \pageref{topologies}). The lattices of admissible subobjects of a
       topological space is the Boolean algebra of subspaces of the
       space, and it transpires $\Int{\Top}{\Top} =
       \pmb{BiTop}$, the category of bitopological spaces and
       bicontinuous maps.

 \item \Loc satisfy the conditions (see \S \ref{locales}, page
       \pageref{locales}). As observed earlier, \Loc is signifivantly
       different from \Set or \Top but there are interesting internal
       topologies on a locale, as investigated in the papers
       \cite{DubeIghedo2016} \& \cite{DubeIghedo2016a}. 

 \item A category \Bb{A} is said to be {\em regular} if it has finite
       limits, kernel pairs have coequalisers and regular epimorphisms
       are pullback stable. Every regular category has a 
       \fact{RegEpi}{Mon} system. Hence every subobject of an object is
       admissible. \\

       A regular category \Bb{A} is said to be {\em coherent} or a {\em
       pre-logos} (see \cite{FreydScedrov1990}) if for each object $X$
       the semilattice \Sub{X}{} is a lattice and for every morphism $f$
       the preimage \finv{f}{} is a lattice homomorphism.\\

       A coherent category \Bb{A} is a {\em Heyting category} or a {\em
       logos} (see \cite{FreydScedrov1990}) or a {\em quasi-category}
       (see \cite{Joyal2008}), if further for every morphism $f$ the
       preimage \finv{f}{} preserve arbitrary joins. It is well known
       from \cite{FreydScedrov1990}, as well as shown in Corollary
       \ref{ppj>subobjectframes} (see page
       \pageref{ppj>subobjectframes}), if each preimage preserve
       arbitrary joins then each lattice of admissible subobjects is a
       frame.\\ 

       In particular, every Heyting category satisfies the
       conditions. Since a topos is an example of a Heyting category,
       every topos satisfies the conditions.

 \item A category \Bb{A} is said to be {\em extensive} if it has finite
       sums and for objects $X$ and $Y$ of \Bb{A}, the canonical functor
       \Arr{+}{\slice{\Bb{A}}{X} \times
       \slice{\Bb{A}}{Y}}{\slice{\Bb{A}}{(X + Y)}} is an equivalence of
       categories (see \cite[Definition
       2.1]{CarboniLackWalters1993}). If further \Bb{A} is small
       complete and small cocomplete then it has an \fact{Epi}{ExtMon}
       system and for each object $X$ of \Bb{A} the lattice
       \Sub{X}{\mathsf{M}} is a distributive complete lattice.\\

       Any quasitopos with disjoint coproducts is extensive; if further
       it has a proper factorisation system then it satisfies the
       conditions.\\

       The categories \Cat of small categories, \opp{\CRing} of affine
       schemes, and the category $\pmb{Sch}$ of schemes are all
       infinitary lextensive with proper factorisation structures. Hence
       they satisfy the conditions. 
       
 \item If \Bb{A} has an \fact{\mathsf{E}}{\mathsf{M}} system then for
       any object $X$ of \Bb{A} the category \slice{\Bb{A}}{X} of
       bundles over $X$ has \fact{\mathsf{E}_X}{\mathsf{M}_X} system,
       where: 
       \begin{align*}
	\mathsf{E}_X = \bigl\{\Arr{e}{(X, x)}{(Y, y)}: e \in
	\mathsf{E}\bigr\} \\ \intertext{and}
	\mathsf{M}_X = \bigl\{\Arr{m}{(X, x)}{(Y, y)}: m \in
	\mathsf{M}\bigr\}.	
       \end{align*}

       If the \fact{\mathsf{E}}{\mathsf{M}} is proper then so also is
       the \fact{\mathsf{E}_X}{\mathsf{M}_X} (see \cite[\S
       2.10]{MMCEGWT2004} for details).\\

       Hence, if \Bb{A} satisfy the conditions of this paper then so
       does each \slice{\Bb{A}}{X}. 
\end{enumerate}

Finally, the regular epimorphisms of internal neighbourhood spaces have
been established in \S \ref{RegEpi} (pages \pageref{RegEpi} -
\pageref{Examples}). Theorem \ref{regepipreobj-desc} (see page
\pageref{regepipreobj-desc}) describes the regular epimorphisms of
internal preneighbourhood spaces. This is similar to the
description of regular epimorphisms of pretopological spaces (see
\cite[Theorem 29, page 16]{BentleyHerrlichLowen1991} and also in
\cite{Kent1969}). The dissimilarity is a consequence of the pullback
stability of epimorphisms of \Set, which is not the case in general (see
\ref{extracondition} \& \ref{pullbackstability>simpledescregepi}, page
\pageref{extracondition}). The pullback stability of epimorphisms in
\Set is also responsible for the regular epimorphisms of pretopological
spaces to be hereditary (see \cite[Theorem 30, page
16]{BentleyHerrlichLowen1991} for the hereditary property). In Theorem
\ref{hereditaryregepidesc} 
(see page \pageref{hereditaryregepidesc}) the hereditary regular 
epimorphisms of internal preneighbourhood spaces are described.

The pullback stability of epimorphisms is a weak condition ensuring
regular epimorphism of internal preneighbourhood spaces to be
hereditary. Theorem \ref{condensurherregepi} (see page
\pageref{condensurherregepi}) provide five conditions which ensure
heredity of regular epimorphisms of internal preneighbourhood
spaces. Regular epimorphisms of preneighbourhood spaces are not stable
under pullbacks --- for instance, regular epimorphisms of pretopological
spaces are not closed under products (see \cite[Example 4, page
5]{BentleyHerrlichLowen1991} for details).

The regular epimorphisms of internal neighbourhood spaces are
similar to the regular epimorphisms of internal preneighbourhood spaces
(see Theorem \ref{regepinbdobj-desc}, page
\pageref{regepinbdobj-desc}). The dissimilarity lies in the replacement
of \Bb{A} by \PPJ{\Bb{A}} over which \NHD{\Bb{A}} is topological. The
hereditary regular epimorphisms of internal neighbourhood spaces is a
little more intricate. Theorem \ref{hereditaryregepiinnbd} (page
\pageref{hereditaryregepiinnbd}) provides alternative characterisations
for hereditary regular epimorphisms of internal neighbourhood
spaces. In the special case of Theorem \ref{hereditaryregepiinnbd}, a
morphism of internal neighbourhood spaces is a regular epimorphism of
internal preneighbourhood spaces  if and only if  it is {\em
pseudo-open} (see \ref{pseudoopendef}, page
\pageref{pseudoopendef} and also \cite[Proposition 28, page
15]{BentleyHerrlichLowen1991} for comparison with \Set).

A summary of the above mentioned connections between different
categories of internal neighbourhood spaces appear in Figure
\ref{catstate} (page \pageref{catstate}). 

The notation and terminology follows \cite{Maclane1997} on
categories and \cite{PicadoPultr2012} on frames and ordered algebraic
systems. 
\section{Preliminaries}
\label{Preliminaries}

\subsection{Factorisation Systems and Admissible Subobjects}
\label{Factorisations}

The modern notion of a factorisation system
\opair{\mathsf{E}}{\mathsf{M}} was introduced in \cite{FreydKelly1972}.
The earlier {\em bicategorical structures} of Maclane (see
\cite[\S 9]{Maclane1950}) can be seen today as those
factorisation systems where 
every $\mathsf{E}$ is an epimorphism and each $\mathsf{M}$ is a
monomorphism. In this part the necessary facts for factorisation
systems are collected from \cite[\S 2]{Janel1997b}.

\subsubsection{Prefactorisation Systems}

Given the morphisms $p$ and $i$ of a category \Bb{A}, the symbol {\em
 \down{p}{i}} is used to denote the statement: if $\comp{v}{p} =
 \comp{i}{u}$, then there exists a unique {\em diagonal
 morphism} $w$ such that the diagram $\xymatrix{ {\cdot} \ar@{->}[r]^{p} 
 \ar@{->}[d]_{u} & {\cdot} \ar@{->}[d]^{v} \ar@{..>}[dl]|{w}\\   
 {\cdot} \ar@{->}[r]_{i} & {\cdot} }$ commutes.

 Let for any set \CAL{H} of morphisms of \Bb{A}:
\begin{equation*}
   \up{\CAL{H}} = \bigl\{p: h \in \CAL{H} \Rightarrow \down{p}{h} \}
    \quad\text{ and }\quad \dn{\CAL{H}} = \bigl\{p: h \in \CAL{H}
    \Rightarrow \down{h}{p}\bigr\}.
\end{equation*}  
 \index{\down{p}{i}}

\begin{Df}
 \label{pfs-df}
 A {\em prefactorisation system} for a category \Bb{A} is 
 a pair \opair{\mathsf{E}}{\mathsf{M}} of sets of morphisms of \Bb{A} 
 such that $\mathsf{E} = \up{\mathsf{M}}$ and $\mathsf{M} =
 \dn{\mathsf{E}}$.
\end{Df}

\begin{Thm}[see {\cite[\S 2]{Janel1997b}}]
 \label{pfs+rel-prop}
 In any category \Bb{A}:
 
 \begin{enumerate}
  \item \label{r+t:up>l-iso} If $f = \comp{m}{e}$, \down{f}{m} and
	\down{e}{m} then $m$ is an isomorphism.
  \item \label{pfs-prop} Given any prefactorisation system
	\opair{\mathsf{E}}{\mathsf{M}} of \Bb{A}:
	\begin{enumerate} 
	 \item \label{iso-common} $\mathsf{E} \cap \mathsf{M} = 
	       \Iso{\Bb{A}}$.   
	 \item \label{compclosed} $\mathsf{M}$ is closed under
	       compositions. 
	 \item \label{M=likemono} If $\comp{g}{f} \in \mathsf{M}$ and
	       $g$ is either a monomorphism or $g \in \mathsf{M}$ then
	       $f \in \mathsf{M}$.
	 \item \label{M=pulbstable} $\mathsf{M}$ is stable under
	       pullbacks. 
	 \item \label{M=limitclosed} $\mathsf{M}$ is closed under
	       limits, i.e., if \Nattrn{\alpha}{F}{G}{\Bb{Z}}{\Bb{A}}
	       with each component of $\alpha$ in $\mathsf{M}$ and both
	       the limits \limit{F}{} and \limit{G}{} exist, then
	       $\limit{\alpha}{} \in \mathsf{M}$.  
	\end{enumerate} 
 \end{enumerate}
\end{Thm}

\subsubsection{Factorisation Systems}

\begin{Df}
 \label{factorisationsys-df}
 A {\em factorisation system} for a category \Bb{A} is a
 prefactorisation system \opair{\mathsf{E}}{\mathsf{M}}, such that
 any morphism $f$ of \Bb{A}, $f = m\circ e$, for some $m \in \mathsf{M}$
 and $e \in \mathsf{E}$.

 An \fact{\mathsf{E}}{\mathsf{M}} for a category \Bb{A} is {\em proper} 
 if $\mathsf{E} \subseteq \Epi{\Bb{A}}$ and $\mathsf{M} \subseteq
 \Mono{\Bb{A}}$.  
\end{Df} 

Every finitely complete and finitely cocomplete category with all
intersections admit a \fact{\Epi{\Bb{A}}}{\ExtMon{\Bb{A}}}.

Furthermore, in any category \Bb{A} with binary products and coproducts
the condition $\mathsf{E} \subseteq \Epi{\Bb{A}}$ implies
$\ExtMon{\Bb{A}} \subseteq \mathsf{M}$ and dually $\mathsf{M} \subseteq
\Mono{\Bb{A}}$ implies $\ExtEpi{\Bb{A}} \subseteq \mathsf{E}$. Hence
such a category  \Bb{A} has a proper \fact{\mathsf{E}}{\mathsf{M}}
implies $\ExtEpi{\Bb{A}}
\subseteq \mathsf{E} \subseteq \Epi{\Bb{A}}$ and $\ExtMon{\Bb{A}}
\subseteq \mathsf{M} \subseteq \Mono{\Bb{A}}$. 

\subsubsection{Admissible Subobjects}

Let \Bb{A} be a finitely complete category with coproducts and a proper 
factorisation system \opair{\mathsf{E}}{\mathsf{M}}.

Given any object $X$ of \Bb{A}, the (possibly large) set
$\bigl\{\Arr{m}{M}{X}: m \in \mathsf{M}\bigr\}$ is endowed with a
natural preorder: 
\[
 m \leq n \Leftrightarrow (\exists p)(m = n\circ p).
\]
The corresponding quotient set is the poset \Sub{X}{\mathsf{M}}
of all {\em $\mathsf{M}$-subobjects} or {\em admissible subobjects} of
$X$. For brevity, an admissible subobject
\Arr{m}{M}{X} of $X$  shall be simply expressed by the morphism
$m$ or even sometimes by $M$.

Since \Bb{A} is finitely complete and $\mathsf{M}$ is pullback
stable, the poset \Sub{X}{\mathsf{M}} of admissible subobjects is
a meet semilattice with largest element. Furthermore the
existence of finite coproducts along with the
\fact{\mathsf{E}}{\mathsf{M}} ensures the existence of finite joins in
\Sub{X}{\mathsf{M}}. Hence \Sub{X}{\mathsf{M}} is a lattice. 

Henceforth in this paper the following stipulation is made on \Bb{A}:

\begin{quotation}
 \label{basicassumption}
 {\em
 \Bb{A} is a finitely complete category with finite coproducts and a
 proper \fact{\mathsf{E}}{\mathsf{M}} such that for each object $X$,
 \Sub{X}{\mathsf{M}} is a distributive complete lattice. 
 }	 
\end{quotation}

Such categories certainly exist. \Set is the most familiar example. In
\S \ref{Examples} (page \pageref{Examples}) it is shown that \Top, \Loc
are examples, and in \S \ref{Introduction} (page 
\pageref{contextlist}) other examples are described.

\subsubsection{Images and Preimages}

From our assumption, for each object $X$ of \Bb{A}, \Sub{X}{\mathsf{M}} is a
distributive complete lattice. The smallest admissible subobject in 
\Sub{X}{\mathsf{M}} is \Arr{\eulf{z}_X}{\emptyset_X}{X} and the largest
is obviously \id{X}. If \Bb{A} has a strict initial object $\emptyset$
then for each object $X$, $\emptyset_X = \emptyset$ --- a situation
familiarly seen in ‌\Set, \Top, \Loc or in any extensive category (see
\cite{CarboniLackWalters1993}). 

Given a morphism \Arr{f}{X}{Y} of \Bb{A} and an admissible subobject
\Arr{n}{N}{Y} of $Y$, the pullback $\xymatrix{ {\finv{f}{N}}
\ar[r]^(0.6){f_n} \ar[d]_{\finv{f}{n}} &  {N} \ar[d]^n \\ {X} \ar[r]_f &
{Y}  }$ of $f$ along $n$ exists. Since $\mathsf{M}$ is pullback
stable (see Theorem
\ref{pfs+rel-prop}\eqref{M=pulbstable}, page \pageref{pfs+rel-prop}),
$\finv{f}{n} \in \Sub{X}{\mathsf{M}}$, and is called the {\em preimage
of $n$ under $f$}. The morphism \Arr{f_n}{\finv{f}{N}}{N} shall be called
the {\em restriction of $f$} to the admissible subobject $N$.  

Given a morphism \Arr{f}{X}{Y} of \Bb{A} and an admissible subobject
\Arr{m}{M}{X} of $X$ the \fact{\mathsf{E}}{\mathsf{M}}:
\[
 \xymatrix@!=12ex{
 {M} \ar[r]^(0.42){\rest{f}{M} (\in \mathsf{E})} \ar[d]_m & {\img{f}{M}}
 \ar[d]^{\img{f}{m} (\in \mathsf{M})} \\
 {X} \ar[r]_f & {Y}
 }
\]
of \comp{f}{m} yields the admissible subobject \img{f}{m} of $Y$, called
the {\em image of $m$ under $f$}. The morphism \rest{f}{M} shall be
called the {\em trace} of $f$ on $M$.

\begin{Thm}
 \label{image-|preimage}
 Given any morphism \Arr{f}{X}{Y} in \Bb{A}, the image
 \Arr{\img{f}{}}{\Sub{X}{\mathsf{M}}}{\Sub{Y}{\mathsf{M}}} and preimage
 \Arr{\finv{f}{}}{\Sub{Y}{\mathsf{M}}}{\Sub{X}{\mathsf{M}}} are order
 preserving maps between the distributive complete lattices of
 admissible subobjects with \adjt{\img{f}{}}{\finv{f}{}}.
\end{Thm}

\begin{Cor}
 \label{image-preimage-consequences}
 Given any morphism \Arr{f}{X}{Y} of \Bb{A}, we have:
 \begin{enumerate}[label=(\alph*),ref=\alph*,align=left]
  \item For every admissible subobject $M$ of $X$, $M \subseteq 
	\finv{f}{\img{f}{M}}$.
  \item For every admissible subobject $N$ of $Y$, $\img{f}{\finv{f}{N}}
	\subseteq N$. 
  \item \img{f}{} preserve all joins and \finv{f}{} preserve all meets.
  \item \label{everyrestepi>imagesurjective} For any admissible
	subobject \Arr{n}{N}{Y} of $Y$, $f_n \in \mathsf{E}$, if and
	only if, $\img{f}{\finv{f}{n}} = n$.
  \item \label{finM>imageinjective} If $f \in \mathsf{M}$ then
	$\finv{f}{\img{f}{M}} = M$, for every admissible subobject $M$
	of $X$.
  \item \label{invwhole=whole+intsurjectivity} $\finv{f}{Y} = X$,
	$\img{f}{\emptyset_X} = \emptyset_Y$ and $f \in \mathsf{E}$, if
	and only if, $Y = \img{f}{X}$.  
 \end{enumerate}
\end{Cor}
\subsection{Filters on an Object}
\label{Filters}

\subsubsection{The Coherent Frame of Filters}

A {\em filter} on an object $X$ of \Bb{A} is just a filter in the
distributive lattice\footnote{A {\em filter} on a meet semilattice $A$
is a subset $F \subseteq A$ which is up-closed (i.e., $x \geq y \in F
\Rightarrow x \in F$) and closed under finite meets (i.e., $x, y \in F
\Rightarrow x \wedge y \in F$).} \Sub{X}{\mathsf{M}} of admissible
subobjects of $X$. The set of all filters on $X$ is \Fil{X} and is
ordered by set theoretic inclusion. 

\begin{Df}
 \label{aux-df} 
 Let $X$ be a bounded lattice.
 \begin{enumerate}[label=(\roman*),align=left]
  \item An element $p$ in a lattice $X$ is said to be {\em compact}, if
	for every subset $S \subseteq X$, there exists a finite subset
	$T \subseteq S$ such that $p \leq \bigvee T$, whenever $p \leq
	\bigvee S$.
  \item A lattice is {\em compact} if its largest element is
	compact.
  \item A lattice is {\em algebraic} if each of its
	elements is the supremum of compact elements.
  \item A frame is {\em coherent} if it is a compact,
	algebraic lattice in which the set of compact elements is closed
	under finite meets.
 \end{enumerate}
 \end{Df}

 \begin{Thm}[{\cite[Theorem 1.2]{IberkleidMcGovern2009b}}]
 \label{upsetmakefreeframe}
 For any object $X$ of \Bb{A}, \Fil{X} is a coherent frame.
 \end{Thm}

 The compact elements of the frame \Fil{X} are precisely \atleast{x} ($x
 \in \Sub{X}{\mathsf{M}}$), where:
 \begin{equation}
  \label{principalfilters-def}
 \atleast{x} = \bigl\{u \in \Sub{X}{\mathsf{M}}: x \leq u\bigr\},   
 \end{equation}
 is the {\em principal} filter on $X$ containing the admissible subobject
 $x$. Clearly:
 \begin{align}
  \label{atleastpreslattop}
  \atleast{(x \vee y)} = (\atleast{x}) \cap (\atleast{y}) & \quad\text{ and }
  & \atleast{(x \wedge y)} = (\atleast{x}) \vee (\atleast{y}).
 \end{align}
 
\subsubsection{Forward and Inverse Filters}

Given a morphism \Arr{f}{X}{Y} of \Bb{A}, filters $A \in \Fil{X}$,
$B \in \Fil{Y}$, let:
\begin{align}
 \label{imageupset-eq}
 \imgfil{f}{A} = \bigl\{y \in \Sub{Y}{\mathsf{M}}: (\exists x \in
 A)(\img{f}{x} \leq y)\bigr\} = \bigl\{y \in \Sub{Y}{\mathsf{M}}:
 \finv{f}{y} \in A\bigr\}, \\ \intertext{and}
 \label{invupset-eq}
 \invfil{f}{B} = \bigl\{x \in \Sub{X}{\mathsf{M}}: (\exists y \in
 B)(\finv{f}{y} \leq x)\bigr\}.
\end{align}
\index{\imgfil{f}{A}}
\index{\invfil{f}{B}}

The filter \imgfil{f}{A} is the smallest filter in \Fil{Y} which contains
the images \img{f}{a} ($a \in A$) and shall be called the {\em forward
filter of $A$ under $f$}. Similarly, \invfil{f}{B} is the smallest
filter in \Fil{X} which contains the preimages \finv{f}{b} ($b \in B$),
and shall be called the {\em inverse filter of $B$ under $f$}.
 
\begin{Thm}
 \label{changeofbase-for-upsets}
 Given any morphism \Arr{f}{X}{Y} of \Bb{A}, both \invfil{f}{} and
 \imgfil{f}{} preserve principal filters, and there is the Galois
 connection: 
 \begin{equation}
  \label{changeofbase-for-upsets-eq}
   \xymatrix@!=18ex{
   {\Fil{X}} \ar@/_2ex/[r]_{\imgfil{f}{}}
   \ar@/^2ex/@{<-}[r]^{\invfil{f}{}} \ar@{}[r]|{\bot} & {\Fil{Y}}
   }.
 \end{equation}

 If further \Arr{\img{f}{}}{\Sub{X}{\mathsf{M}}}{\Sub{Y}{\mathsf{M}}}
 preserve finite meets then \Arr{\imgfil{f}{}}{\Fil{X}}{\Fil{Y}} has a
 right adjoint
 \Arr{\Imgfil{f}{}}{\Sub{Y}{\mathsf{M}}}{\Sub{X}{\mathsf{M}}} defined
 by:
  \[
 \Imgfil{f}{B} = \bigl\{x \in \Sub{X}{\mathsf{M}}: \img{f}{x}
 \in B\bigr\}.
 \]
\end{Thm}

\begin{proof}
 Given the filters $A \in \Fil{X}$, $B \in \Fil{Y}$:
 \begin{multline*}
  \invfil{f}{B} \subseteq A \Leftrightarrow (\forall b \in
  B)(\finv{f}{b} \in A) \\
  \Leftrightarrow (\forall b \in B)(\exists a \in A)(a \leq
  \finv{f}{b}) \\
  \Leftrightarrow (\forall b \in B)(\exists a \in A)(\img{f}{a} \leq
  b) 
  \Leftrightarrow B \subseteq \imgfil{f}{A},
 \end{multline*}
 proving \adjt{\invfil{f}{}}{\imgfil{f}{}}.

 On the other hand, given any family \seq{A}{i}{I} of filters on $X$ if
 \img{f}{} preserve finite meets then:
 \begin{multline*}
  b \in \imgfil{f}{\bigl(\bigvee_{i \in I} A_i\bigr)} \Leftrightarrow \finv{f}{b} \in
  \bigvee_{i \in I}A_i \\
  \Leftrightarrow (\exists n \geq 1)(\exists i_1, i_2, \dots, i_n \in
  I)(\exists a_1 \in A_{i_{1}}, a_2 \in A_{i_{2}}, \dots a_n \in
  A_{i_{n}})\bigl(a_1 \wedge a_2 \dots \wedge a_n \leq \finv{f}{b}\bigr) \\
  \Leftrightarrow (\exists n \geq 1)(\exists i_1, i_2, \dots, i_n \in
  I)(\exists a_1 \in A_{i_{1}}, a_2 \in A_{i_{2}}, \dots a_n \in
  A_{i_{n}})\bigl(\img{f}{(a_1 \wedge a_2 \dots \wedge a_n)} 
  \leq b\bigr) \\
  \Leftrightarrow (\exists n \geq 1)(\exists i_1, i_2, \dots, i_n \in
  I)(\exists a_1 \in A_{i_{1}}, a_2 \in A_{i_{2}}, \dots a_n \in
  A_{i_{n}})\bigl(\img{f}{a_1} \wedge \img{f}{a_2} \wedge \dots \wedge \img{f}{a_n} 
  \leq b\bigr) \\
  \Leftrightarrow b \in \bigvee_{i \in I}\imgfil{f}{A_i},
 \end{multline*}
 indicate \imgfil{f}{} preserve all joins and hence must have
 a right adjoint \Arr{\Imgfil{f}{}}{\Fil{Y}}{\Fil{X}}.

 Finally, for any admissible subobject $x$ of $X$:
 \begin{equation*}
  x \in \Imgfil{f}{B}
 \Leftrightarrow \atleast{x} \subseteq \Imgfil{f}{B} \Leftrightarrow
 \imgfil{f}{(\atleast{x})} \subseteq B \Leftrightarrow
 \atleast{(\img{f}{x})} \subseteq B \Leftrightarrow \img{f}{x} \in B
 \end{equation*}
 implying:   
 \[
 \Imgfil{f}{B} = \bigl\{x \in \Sub{X}{\mathsf{M}}: \img{f}{x}
 \in B\bigr\}.
 \]

\end{proof}

\subsection{Preimage preserving joins}

Sometimes the preimage
\Arr{\finv{f}{}}{\Sub{Y}{\mathsf{M}}}{\Sub{X}{\mathsf{M}}} for a
morphism \Arr{f}{X}{Y} of \Bb{A} preserve joins --- for instance in
\Set, \Top, \Meas, and in many other concrete categories. However, in
\Loc the preimages usually preserve finite joins only (see
\cite[Proposition 9.2, page 41]{PicadoPultr2012}).

\begin{Thm}
 \label{ppj full or partial}
 \begin{enumerate}[label=(\alph*),ref=(\alph*),align=left]
  \item \label{monoppj<>meetdistoverjoins}
	If \Arr{p}{P}{X} be an admissible subobject of $X$ then the
	preimage function
	\Arr{\finv{p}{}}{\Sub{X}{\mathsf{M}}}{\Sub{P}{\mathsf{M}}} 
	preserve arbitrary joins, if and only if, for every family
	\seq{m}{i}{I} of admissible subobjects of $X$:
	\begin{equation}
	 \label{frame law}
	  p \wedge \bigvee_{i \in I}m_i = \bigvee_{i \in I}(p \wedge m_i).
	\end{equation}

  \item  \label{preimagehasrightadjoint-alt}
	 The following are equivalent for any morphism \Arr{f}{X}{Y} of
	 \Bb{A}: 
	 \begin{enumerate}[label=(\roman*),ref=(\roman*),align=left] 
	  \item \label{preimagehasrightadjoint}
		\Arr{\finv{f}{}}{\Sub{Y}{\mathsf{M}}}{\Sub{X}{\mathsf{M}}}
		has a right adjoint
		\Arr{\fInv{f}{}}{\Sub{X}{\mathsf{M}}}{\Sub{Y}{\mathsf{M}}}.   
	  \item \label{preimagepresarbjoins}
		\Arr{\finv{f}{}}{\Sub{Y}{\mathsf{M}}}{\Sub{X}{\mathsf{M}}}
		preserve all arbitrary joins. 
	  \item \label{invfilpreservearbmeets}
		\Arr{\invfil{f}{}}{\Fil{Y}}{\Fil{X}} preserve all
		arbitrary meets. 	 
	  \item \label{invfilhasleftadjoint}
		\Arr{\invfil{f}{}}{\Fil{Y}}{\Fil{X}} has a left adjoint
		\Arr{\Invfil{f}{}}{\Fil{X}}{\Fil{Y}}. 
	 \end{enumerate}
 \end{enumerate}
\end{Thm}

\begin{proof}
 \begin{enumerate}[label=(\alph*),align=left]
  \item If \Arr{\finv{p}{}}{\Sub{X}{\mathsf{M}}}{\Sub{P}{\mathsf{M}}}
	preserve arbitrary joins, then for any family \seq{m}{i}{I} of
	admissible subobjects of $X$:
	\begin{multline*}
	 p \wedge \bigvee_{i \in I}m_i = p\circ\finv{p}{\bigl(\bigvee_{i
	 \in I}m_i\bigr)} \\
	 = p\circ\bigl(\bigvee_{i \in I}\finv{p}{m_i}\bigr) 
	 = \img{p}{\bigl(\bigvee_{i \in I}\finv{p}{m_i}\bigr)} \\
	 = \bigvee_{i \in I}\img{p}{\finv{p}{m_i}} = \bigvee_{i \in 
	 I}(p\circ \finv{p}{m_i}) = \bigvee_{i \in I}(p \wedge m_i),
	\end{multline*}
	proving \eqref{frame law}.
	
	On the other hand, if \eqref{frame law} holds, then:
	\begin{multline*}
	 p\circ\finv{p}{\bigl(\bigvee_{i \in I}m_i\bigr)} = p \wedge
	 \bigvee_{i \in I}m_i = \bigvee_{i \in I}(p \wedge m_i) \\
	 = \bigvee_{i \in I}\img{p}{\finv{p}{m_i}} =
	 \img{p}{\bigl(\bigvee_{i \in I}\finv{p}{m_i}\bigr)} =
	 p\circ\bigl(\bigvee_{i \in I}\finv{p}{m_i}\bigr),
	\end{multline*}
	which implies $\finv{p}{\bigl(\bigvee_{i \in I}m_i\bigr)} =
	\bigvee_{i \in I}\finv{p}{m_i}$, completing the proof for this
	part. 

  \item  Obviously, \ref{preimagehasrightadjoint} and
	 \ref{preimagepresarbjoins} are equivalent and so also for the
	 pair \ref{invfilpreservearbmeets} and
	 \ref{invfilhasleftadjoint}.  

	 Assuming \ref{preimagepresarbjoins}, given a family
	 \seq{B}{i}{I} of filters on $Y$, $p \in \bigcap_{i \in
	 I}\invfil{f}{B_i}$, if and only if, for each $i \in I$ there
	 exist a $b_i \in B_i$ such that $\finv{f}{b_i} \leq p$, which
	 imply $\bigvee_{i \in I}\finv{f}{b_i} =
	 \finv{f}{\bigl(\bigvee_{i \in I}b_i\bigr)} \leq p \Rightarrow p
	 \in \invfil{f}{\bigl(\bigcap_{i \in I}B_i\bigr)}$, since
	 $\bigvee_{i \in I}b_i \in \bigcap_{i \in I}B_i$, showing
	 \ref{invfilpreservearbmeets} follows.  

	 On the other hand, assuming \ref{invfilpreservearbmeets} and
	 using it on principal filters shows \ref{preimagepresarbjoins}
	 to follow. 
 \end{enumerate}
 \end{proof}

\begin{enumerate}[itemsep=1.2ex,label=\underline{Remark
 \arabic*.},ref=\underline{Remark \arabic*.},align=left,resume=rem]
 \item A morphism \Arr{f}{X}{Y} of \Bb{A} for which the preimage \Arr{\finv{f}{}}{\Sub{Y}{\mathsf{M}}}{\Sub{X}{\mathsf{M}}} preserve arbitrary joins shall be said to have
       {\em preimage preserve joins} property. For any such morphism: 
       \begin{align}
	\label{rtadjointofpreimage-eq}
	\fInv{f}{x} = \bigvee\bigl\{y \in \Sub{Y}{\mathsf{M}}:
	\finv{f}{y} \leq x\bigr\}, \quad\text{ for }x \in
	\Sub{X}{\mathsf{M}}, \\ 
	\intertext{and}
	\label{lftadjointofinvfil-eq}
	\Invfil{f}{A} = \bigl\{y \in \Sub{Y}{\mathsf{M}}: (\exists a
	\in A)(\fInv{f}{a} \leq y)\bigr\}, \quad\text{ for }A \in
	\Fil{X}. 
       \end{align}

       In particular, \Invfil{f}{} preserve principal filters and:
       \begin{equation}
	\label{lftadjointofinvfilpresprincipalfilters}
	 \Invfil{f}{\atleast{x}} = \atleast{\fInv{f}{x}}, \quad\text{ for }x
	 \in \Sub{X}{\mathsf{M}}.
       \end{equation}       
\end{enumerate}

As a consequence of Theorem \ref{ppj full or
partial}\ref{monoppj<>meetdistoverjoins}: 

\begin{Cor}
 \label{ppj>subobjectframes}
 If every morphism (or, every admissible monomorphism) of \Bb{A} has the
 preimage preserve join property then each \Sub{X}{\mathsf{M}} is a frame.
\end{Cor}

The results in this section are well known and are also seen in 
\cite{FreydScedrov1990}.     % introduction to \S 1.2 to be given
\section{Preneighbourhoods, Weak Neighbourhoods and Neighbourhoods}
\label{Neighbourhoods}

In this section we shall define the notion of a {\em neighbourhood} of an
admissible subobject and develop some of their relevant properties.

\subsection{Neighbourhoods}
\label{neighbourhoodobjects_v2}

\begin{Df}
 \label{neighbourhoodstructures-df}
 Given an object $X$ of \Bb{A}:
 \begin{enumerate}[itemsep=1.2ex,label=(\alph*),align=left,ref=(\alph*)]
  \item \label{preneighbourhood-df}
	An order preserving map
	\Arr{\mu}{\opp{\Sub{X}{\mathsf{M}}}}{\Fil{X}} is a {\em
	preneighbourhood } on $X$ if for each $m \in
	\Sub{X}{\mathsf{M}}$:   
	\begin{equation}
	 \label{upmajorises}
	  n \in \mu(m) \Rightarrow m \leq n. 
	\end{equation}
	 
  \item \label{weakneighbourhood-df}
	A preneighbourhood  $\mu$ on $X$ is a {\em weak
	neighbourhood } on $X$ if:
	\begin{equation}
	 \label{seqofseqcondition}
	  \mu(m) \subseteq \bigcup_{p \in \mu(m)}\bigcap_{x \leq
	  p}\mu(x), \quad\text{ for }m \in \Sub{X}{\mathsf{M}}.
	\end{equation}
	
  \item \label{neighbourhood-df}
	A weak neighbourhood  $\mu$ on $X$ is a {\em
	neighbourhood } on $X$ if:
	\begin{equation}
	 \label{presarbmeets}
	  \mu(\bigvee G) = \bigcap_{x \in G}\mu(x), \quad\text{ for all
	  }G \subseteq \Sub{X}{\mathsf{M}}.
	\end{equation}
 \end{enumerate}
	\end{Df}

\begin{enumerate}[itemsep=1.2ex,label=\underline{Remark
 \arabic*},ref=\underline{Remark \arabic*},align=left,resume=rem]
 \item The set of all preneighbourhoods, weak neighbourhoods or
       neighbourhoods on and object $X$ is 
       denoted by the symbols \pnhd{X}, \wnhd{X} or \nhd{X},
       respectively.\\

       If $\mu$ is a preneighbourhood (respectively, weak neighbourhood,
       neighbourhood) on and object $X$ of \Bb{A} then the pair
       \opair{X}{\mu} shall be called an {\em internal preneighbourhood
       space} (respectively, {\em internal weak neighbourhood space},
       {\em internal neighbourhood space}).
       
 \item \label{trivialnbds}
       Surely, \Arr{\atleast{}}{\opp{\Sub{X}{\mathsf{M}}}}{\Fil{X}},
       where $\atleast{m} = \bigl\{p \in \Sub{X}{\mathsf{M}}: m \leq 
       p\bigr\}$ for any $m \in \Sub{X}{\mathsf{M}}$, is
       a neighbourhood  on $X$. \\

       Further, \Arr{\nabla}{\opp{\Sub{X}{\mathsf{M}}}}{\Fil{X}} defined
       by:
       \begin{equation}
	\label{indiscrete-eq}
	 \nabla(m) =
	 \begin{cases}
	  \Sub{X}{\mathsf{M}}, & \text{ if }m= \emptyset_X \\
	  \bigl\{\id{X}\bigr\}, & \text{ otherwise}
	 \end{cases},
       \end{equation}
       is also a neighbourhood  on $X$. \\

 \item \label{pointwiseorderonpnbd}
       The set \pnhd{X} is ordered {\em pointwise},
       i.e., given preneighbourhoods  $\mu$ and $\nu$ on $X$,
       $\mu \leq \nu$ if for each $m \in \Sub{X}{\mathsf{M}}$, $\mu(m)
       \subseteq \nu(m)$. \\   

       Consequently, \eqref{upmajorises} equivalently states
       $\nabla \leq \mu \leq \atleast{}$ for every preneighbourhood
        $\mu$ on $X$, i.e., \pnhd{X} is a bounded poset.\\  
\end{enumerate}

\subsubsection{Weak Neighbourhoods are Interpolative Preneighbourhoods}

\begin{Thm}
 \label{pnbdisnbd<>interpolative}
 A preneighbourhood  $\mu$ on an object $X$ of
 \Bb{A} is a weak neighbourhood if and only if  it is {\em
 interpolative}, i.e., the following equation holds:
 \begin{equation}
  \label{interpolative-eq}
   \mu(m) = \bigl\{p \in \Sub{X}{\mathsf{M}}: (\exists q \in \mu(m))(p
   \in \mu(q))\bigr\}, \quad\text{ for all }m \in \Sub{X}{\mathsf{M}}.
 \end{equation}

\end{Thm} 
\index{preneighbourhood!interpolative}

\begin{proof}
 Firstly, if $\mu$ be a
 preneighbourhood  on $X$, then for any admissible subobject
 $p$ of $X$, the set $\bigl\{\mu(x): x \leq p\bigr\}$
 of filters on $X$ has $\mu(p)$ as the smallest filter. Hence:
 \[
  \mu(p) = \bigcap_{x \leq p} \mu(x).
 \]

 On the other hand, if  $p \in \mu(m)$ then using \eqref{upmajorises},
 $\mu(p) \subseteq \mu(m)$, and hence:  
 \[
  \mu(m) \supseteq \bigcup_{p \in \mu(m)} \mu(p).
 \]

 Consequently, $\mu$ is a weak neighbourhood  on $X$, if and only
 if, for each admissible subobject $m$ of $X$:
 \[
  \mu(m) \subseteq \bigcup_{p \in \mu(m)}\bigcap_{x \leq p}\mu(x) =
 \bigcup_{p \in \mu(m)}\mu(p) \subseteq \mu(m) \Longrightarrow \mu(m) =
 \bigcup_{p \in \mu(m)}\mu(p),
 \]
 completing the proof.
\end{proof}

\begin{enumerate}[itemsep=1.2ex,label=\underline{Remark \arabic*},ref=\underline{Remark \arabic*},align=left,resume=rem]
 \item \label{filteredjoinoffilters}
       For any preneighbourhood  $\mu$ on $X$ and any $m \in
       \Sub{X}{\mathsf{M}}$, $\bigcup_{p \in \mu(m)}\mu(p)$
       is a subset of $\mu(m)$. Theorem \ref{pnbdisnbd<>interpolative}
       asserts that a preneighbourhood $\mu$ is a weak neighbourhood, if
       and only if, the subset $\bigcup_{p \in \mu(m)}\mu(p)$ is a
       filter and:
       \[
	\mu(m) = \bigvee_{p \in \mu(m)}\mu(p) = \bigcup_{p \in \mu(m)}\mu(p).
       \]

\end{enumerate}

\subsubsection{Complete Lattices \pnhd{X}, \wnhd{X}}

\begin{Thm}
 \label{nbdstructuresmakecompletelattice}
 The set \pnhd{X} of all preneighbourhoods on $X$ with the
 pointwise order in \ref{pointwiseorderonpnbd} (see page
 \pageref{pointwiseorderonpnbd}) is a complete lattice. 

 The subset \wnhd{X} of all weak neighbourhoods on $X$ is also
 a complete lattice with joins computed as in \pnhd{X}. 
\end{Thm}

\begin{proof}
 As observed in \ref{pointwiseorderonpnbd} on page
 \pageref{pointwiseorderonpnbd}, \pnhd{X} is already a bounded poset.

 Furthermore, given any set $T \subseteq \pnhd{X}$ of preneighbourhoods
 on $X$, let for each $m \in \Sub{X}{\mathsf{M}}$:
 \begin{align}
  \label{pnbdsup-df}
  \bigl(\bigvee T\bigr)(m) = \bigvee_{\tau \in T}\tau(m), \\
  \intertext{and}
  \label{pnbdinf-df}
  \bigl(\bigwedge T\bigr)(m) = \bigcap_{\tau \in T}\tau(m).
 \end{align}

 Clearly $\bigvee T, \bigwedge T \in \pnhd{X}$ and $\bigvee T$
 (respectively, $\bigwedge T$) is the supremum (respectively, infimum)
 of $T$ in \pnhd{X}. Hence \pnhd{X} is a complete lattice.

 Given the subset $T \subseteq \wnhd{X}$, $\bigvee T$ is a
 preneighbourhood on $X$. If $p \in \bigl(\bigvee T\bigr)(m)$, then
 there exists a natural number $n \geq 1$, $\tau_1, \tau_2, \dots,
 \tau_n \in T$, $p_1 \in \tau_1(m)$, $p_2 \in \tau_2(m)$, \dots, $p_n
 \in \tau_n(m)$ such that $p = p_1 \wedge p_2 \wedge \dots \wedge
 p_n$. Using Theorem \ref{pnbdisnbd<>interpolative}, there exist $q_1
 \in \tau_1(m)$, $q_2 \in \tau_2(m)$, \dots, $q_n \in \tau_n(m)$ such
 that $p_1 \in \tau_1(q_1)$, $p_2 \in \tau_2(q_2)$, \dots, $p_n
 \in \tau_n(q_n)$. Then $q = q_1 \wedge q_2 \wedge \dots \wedge q_n \in
 \bigl(\bigvee T\bigr)(m)$ and
 \[
  p = p_1 \wedge p_2 \wedge \dots \wedge p_n \in \tau_1(q_1) \vee
 \tau_2(q_2) \vee \dots \vee \tau_n(q_n) \subseteq \tau_1(q) \vee
 \tau_2(q) \vee \dots \vee \tau_n(q) \subseteq \bigl(\bigvee T\bigr)(q)
 \]
 shows $\bigl(\bigvee T\bigr)$ to be interpolative. Hence by Theorem
 \ref{pnbdisnbd<>interpolative}, $\bigvee T \in \wnhd{X}$ and is the
 supremum of $T$ in \wnhd{X}.

 Since \wnhd{X} is a bounded poset with every subset having a supremum,
 it is a complete lattice.
\end{proof}

\begin{enumerate}[itemsep=1.2ex,label=\underline{Remark \arabic*},ref=\underline{Remark \arabic*},align=left,resume=rem]
 \item For any preneighbourhood $\mu$ on $X$, the largest weak
       neighbourhood on $X$ smaller than $\mu$ is:
       \begin{equation}
	\label{cannbd-df}
	 \cannbd{\mu} = \bigvee\bigl\{\nu \in \wnhd{X}: \nu \leq
	 \mu\bigr\}.
       \end{equation}
 \item For any subset $T \subseteq \wnhd{X}$:
       \begin{equation}
	\label{wnbdinf-df}
	 \bigl(\bigwedge T\bigr) = \cannbd{\bigl(\bigcap T\bigr)}.
       \end{equation}
\end{enumerate}

\subsubsection{Open Subobjects and Interiors}

Given a preneighbourhood  $\mu$ on $X$ it is easy to observe for any $p
\in \Sub{X}{\mathsf{M}}$ the three statements in:  
\begin{align}
 \label{nbditself}
 p \in \mu(p), \\
 \label{nbdofallitcontains}
 (\forall m \in \Sub{X}{\mathsf{M}})\bigl(m \leq p \Rightarrow p \in
 \mu(m)\bigr), \\ \intertext{and}
 \label{largestnbdfilter}
 \mu(p) = \atleast{p}
\end{align}
are equivalent: given \eqref{nbditself}, if $m \leq p$ then $p \in
\mu(p) \subseteq \mu(m)$ proves \eqref{nbdofallitcontains}, given
\eqref{nbdofallitcontains}, if $q \geq p$ then $q \in \mu(p)$ proves
\eqref{largestnbdfilter}, and \eqref{largestnbdfilter} automatically
implies \eqref{nbditself}.

Let:

\begin{multline}
 \label{opensets}
 \mo_\mu = \bigl\{p \in \Sub{X}{\mathsf{M}}: p \in \mu(p)\bigr\} \\
 = \bigl\{p \in \Sub{X}{\mathsf{M}}: m \leq p \Leftrightarrow p \in
 \mu(m)\bigr\} \\
 = \bigl\{p \in \Sub{X}{\mathsf{M}}: \mu(p) = \atleast{p}\bigr\}
\end{multline}

and 
\begin{equation}
 \label{interior-eq}
 \intr{\mu}{m} = \bigvee\bigl\{p \in \mo_\mu: p \leq m\bigr\},
 \quad\text{ for }m \in \Sub{X}{\mathsf{M}}.
\end{equation}

The admissible subobjects in $\mo_\mu$ are called {\em $\mu$-open}
subobjects of $X$; for any admissible
subobject $m \in \Sub{X}{\mathsf{M}}$, the admissible subobject
\intr{\mu}{m} is {\em $\mu$-interior} of $m$.

Observe: for any preneighbourhood $\mu$ on $X$, the  largest admissible
subobject \id{X} is always $\mu$-open, and from \eqref{presarbmeets}
(page \pageref{presarbmeets}), if $\mu$ is a neighbourhood on  $X$ then
the smallest subobject $\emptyset_X$ is $\mu$-open. Using
\eqref{atleastpreslattop} (page \pageref{atleastpreslattop}) it follows
that the set $\mo_\mu$ of $\mu$-open subobjects closed under finite
meets. Furthermore, for any preneighbourhood $\mu$ on $X$,
\Arr{\intr{\mu}{}}{\Sub{X}{\mathsf{M}}}{\Sub{X}{\mathsf{M}}} is an order 
preserving idempotent function fixing every $\mu$-open subobjects such
that $\intr{\mu}{m} \leq m$ ($m \in \Sub{X}{\mathsf{M}}$).

\begin{Thm}
 \label{interiorisopen-alt}
 If \Arr{\mu}{\opp{\Sub{X}{\mathsf{M}}}}{\Fil{X}} be a
 preneighbourhood  on $X$ then the set $\mo_\mu$ of $\mu$-open
 subobjects is closed under arbitrary joins if and only if  for every
 $m \in \Sub{X}{\mathsf{M}}$ its $\mu$-interior $\intr{\mu}m$ is
 $\mu$-open.

 Furthermore, in such a case the following two statements are
 equivalent:
 \begin{enumerate}[itemsep=1ex,label=(\alph*),align=left]
  \item \label{opennbdsgenerate}
	For any $m \in \Sub{X}{\mathsf{M}}$, $\mu(m) =
	\bigcup\bigl\{\atleast{q}: m \leq q \in \mo_\mu\bigr\}$.  
  \item \label{nbdiffintcontain} 
	For any $m \in \Sub{X}{\mathsf{M}}$, $p \in \mu(m)
	\Leftrightarrow m \leq \intr{\mu}{p}$. 
 \end{enumerate} 
 
\end{Thm}

\begin{proof}
 The \underline{only if} part of the first statement is immediate from
 the definition of $\mu$-interior in \eqref{interior-eq}. 

 Conversely, if for each $m \in \Sub{X}{\mathsf{M}}$, the $\mu$-interior
 \intr{\mu}{m} of $m$ is $\mu$-open then for any $T \subseteq \mo_\mu$,
 $\intr{\mu}{\bigl(\bigvee T\bigr)} \in \mo_\mu$.   

 Since the elements of $\mo_\mu$ are fixed points of $\mu$-interior
 assignment:
\begin{multline*}
 t \in T \Rightarrow t \leq \bigvee T 
 \Rightarrow t = \intr{\mu}{t} \leq \intr{\mu}{\bigvee T} \\
 \Rightarrow \bigvee T \leq \intr{\mu}{\bigvee T} \Rightarrow \bigvee T
 = \intr{\mu}{\bigvee T} \in \mo_\mu.
\end{multline*}

 For the second part of the statement, assume $\mu$ is a
 preneighbourhood on $X$ such that every $\mu$-interior is
 $\mu$-open. In this case, using \eqref{opensets}:
 \begin{equation*}
  m \leq \intr{\mu}{p} \leq p \Rightarrow p \in \mu(\intr{\mu}{p})
  \subseteq \mu(m),
 \end{equation*}
 implies the \underline{$\Leftarrow$} part of the statement in
 \ref{nbdiffintcontain} is true.

 The implication of \ref{opennbdsgenerate} from \ref{nbdiffintcontain}
 is trivial.
 
 Assuming \ref{opennbdsgenerate}:
 \[
  p \in \mu(m) \Leftrightarrow (\exists q \in \mo_\mu)(m \leq q \leq p)
 \Rightarrow m \leq \intr{\mu}{p},
 \]
 \ref{nbdiffintcontain} follows, completing the proof. 
\end{proof} 
 
\begin{enumerate}[itemsep=1.2ex,label=\underline{Remark \arabic*},ref=\underline{Remark \arabic*},align=left,resume=rem]
 \item Observe, for any preneighbourhood $\mu$ on $X$, since
       $\mu$-interior have $\mu$-open subobjects as fixed points: 
       \begin{equation}
	\label{intislargestopeninside}
	 (\forall p \in \mo_\mu)\bigl(p \leq m \Leftrightarrow p \leq
	 \intr{\mu}{m}\bigr). 
       \end{equation} 

       In the special case when every $\mu$-interior is $\mu$-open, 
       \eqref{intislargestopeninside} provides the familiar meaning:
       \intr{\mu}{m} is the largest $\mu$-open subobject contained in
       $m$.

 \item If for a preneighbourhood $\mu$ every $\mu$-interior is
       $\mu$-open then for admissible subobjects $m$ and $n$ of an
       object $X$, $m \wedge n \geq (\intr{\mu}{m} \wedge \intr{\mu}{n})
       \in \mo_\mu \Rightarrow \intr{\mu}{(m \wedge n)} \geq
       (\intr{\mu}{m} \wedge \intr{\mu}{n})$. \\

       The order preserving property of interior returns:
       \[
	\intr{\mu}{(m \wedge n)} = \intr{\mu}{m} \wedge \intr{\mu}{n}.
       \]
\end{enumerate}

This leads to:

\begin{Cor}
 \label{intropen>Kuratowski}
 The interior operation \intr{\mu}{} of a preneighbourhood $\mu$
 for which the $\mu$-interiors are $\mu$-open is a Kuratowski
 interior operation.
\end{Cor}

\begin{enumerate}[itemsep=1.2ex,label=\underline{Remark \arabic*},ref=\underline{Remark \arabic*},align=left,resume=rem]
 \item Any preneighbourhood  on $X$ which satisfies the condition in
       Theorem \ref{interiorisopen-alt}\ref{opennbdsgenerate} is 
       interpolative and hence a weak neighbourhood .\\

 \item The preneighbourhoods
       \Arr{\mu}{\opp{\Sub{X}{\mathsf{M}}}}{\Fil{X}} satisfying the
       condition in Theorem \ref{interiorisopen-alt}\ref{opennbdsgenerate} are
       a very special kind of weak neighbourhoods --- the ones
       which are determined completely by the $\mu$-open
       subobjects. Using Theorem \ref{snbdgenbyopens} these are in
       between the weak neighbourhoods and neighbourhoods.         
\end{enumerate}

\subsubsection{Interiors for Neighbourhoods}

\begin{Thm}
 \label{snbdgenbyopens}
 If \Arr{\mu}{\opp{\Sub{X}{\mathsf{M}}}}{\Fil{X}} be a neighbourhood
  on $X$ then every $\mu$-interior is $\mu$-open and
 \[
  \mu(m) = \bigcup\bigl\{\atleast{q}: m \leq q \in \mo_\mu\bigr\}.
 \]
\end{Thm}

\begin{proof}
 If $T \subseteq \mo_\mu$ then $\mu(\bigvee T) = \bigcap_{t \in T}\mu(t)
 = \bigcap_{t \in T}\atleast{t} = \atleast{\bigl(\bigvee T\bigr)}$,
 shows the set $\mo_\mu$ of all $\mu$-open subobjects closed
 under arbitrary joins. Using Theorem \ref{interiorisopen-alt} (page
 \pageref{interiorisopen-alt}) completes the proof of the first
 statement.

 Choose and fix a $p \in \mu(m)$. Let: 
 \[
 T_p = \bigl\{u \in \Sub{X}{M}: p \in \mu(u)\bigr\}\quad\text{
 and }p_0 = \bigvee T_p.
 \]
 
 By definition of $T_p$, $m \in T_p$. Since $\mu$ is a
 neighbourhood:
 \[
  \mu(p_0) = \bigcap_{u \in T_p}\mu(u) \Rightarrow p \in \mu(p_0)
 \Rightarrow p_0 \in T_p.
 \]

 Since $\mu$ is a weak neighbourhood 
 on $X$, from Theorem \ref{pnbdisnbd<>interpolative} (page
 \pageref{pnbdisnbd<>interpolative}) it is interpolative. Hence,
 if $u \in T_p$, there exists a $v \in
 \mu(u)$ such that $p \in \mu(v)$. Consequently, $v \in T_p$ and
 $p_0 \geq v \in \mu(u) \Rightarrow p_0 \in \mu(u)$. Thus $p_0
 \in \mu(p_0)$, yielding:
 \[
 p \in \mu(m) \Longrightarrow (\exists q \in \mo_\mu)(m \leq q
 \leq p).
 \]  
 
 Since $\mu$ is a weak neighbourhood:
 \[
 \mu(m) \subseteq \bigcup\bigl\{\atleast{q}: m \leq q \in
 \mo_\mu\bigr\} = \bigcup\bigl\{\mu(q): m \leq q \in
 \mo_\mu\bigr\} \subseteq \bigcup\bigl\{\mu(t): t \in
 \mu(m)\bigr\} = \mu(m)
 \]
 completing the proof.  
\end{proof}

\subsubsection{Kuratowski Interiors and Neighbourhoods}

Let \low{\CAL{K}}{X} be the set of all Kuratowski operations on an
object $X$. This set can be ordered pointwise, producing a partially
ordered set. Using \ref{trivialnbds} (see page \pageref{trivialnbds}),
Theorem \ref{snbdgenbyopens} (page \pageref{snbdgenbyopens}) and
Corollary in \ref{intropen>Kuratowski} (page
\pageref{intropen>Kuratowski}), $\intr{\nabla}{}, \intr{\atleast{}}{}
\in \CAL{K}_X$. 

Let \low{\CAL{P}}{X} be the subset of all preneighbourhoods $\mu$ on $X$
for which the $\mu$-interiors are $\mu$-open. The order from \pnhd{X}
restricts to provide another partially ordered set.

\begin{Thm}
 \label{Kuratowskiisretractofpnbdwithintopen}
 The interior operation
 \Arr{\intr{}{}}{\low{\CAL{P}}{X}}{\low{\CAL{K}}{X}} is a split 
 epimorphism of bounded posets having a left adjoint which restricts to
 an isomorphism precisely on \nhd{X}.  

 Furthermore, for each Kuratowski interior operation $i$ on $X$ the
 fibre \finv{\intr{}{}}{i} has exactly one neighbourhood on $X$,
 which is the smallest element of the fibre.
\end{Thm}

\begin{proof}
 Since $\mo_\nabla = \bigl\{\emptyset_X, X\bigr\}$ (see
 \eqref{indiscrete-eq}, page \pageref{indiscrete-eq}) and
 $\mathfrak{O}_{\atleast{}} = \Sub{X}{\mathsf{M}}$ one obtains:  
 \begin{align*}
  \intr{\nabla}{m} = \begin{cases}
		      \emptyset_X, & \text{ if } m \neq \id{X} \\
		      X, & \text{ otherwise} 
		     \end{cases}, & \quad \text{ and } \quad &
  \intr{\atleast{}}{m} = m,
 \end{align*}
 it follows that \intr{}{} preserve the bounds.

 If $\mu, \nu \in \low{\CAL{P}}{X}$ with $\mu \leq \nu$ then $p \in \mo_\mu
 \Leftrightarrow p \in \mu(p) \subseteq \nu(p) \Rightarrow p \in
 \mo_\nu$ implying $\mo_\mu \subseteq \mo_\nu$.

 Hence $\intr{\mu}{m} = \bigvee\bigl\{p \in \mo_\mu: p \leq
 m\bigr\} \leq \bigvee\bigl\{p \in \mo_\nu: p \leq m\bigr\} =
 \intr{\nu}{m}$, showing
 \Arr{\intr{}{}}{\low{\CAL{P}}{X}}{\low{\CAL{K}}{X}} to be a morphism of
 bounded posets.

 Let \Arr{i}{\Sub{X}{\mathsf{M}}}{\Sub{X}{\mathsf{M}}} be a Kuratowski
 interior operation on $X$ and let:
 \begin{equation}
  \label{kuratowski>nbd}
   \mathfrak{p}_i(m) =  \bigl\{p \in \Sub{X}{\mathsf{M}}: m \leq
   i(p)\bigr\}, \quad\text{ for }m \in \Sub{X}{\mathsf{M}}.
 \end{equation}

 Clearly: $\mathfrak{p}_i(\emptyset_X) = \Sub{X}{\mathsf{M}}$,
 $\mathfrak{p}_i(\id{X}) = \bigl\{\id{X}\bigr\}$, $m \leq n \Rightarrow 
 \mathfrak{p}_i(n) \subseteq \mathfrak{p}_i(m)$, $q \geq p \in
 \mathfrak{p}_i(m) \Rightarrow q \in  
 \mathfrak{p}_i(m)$ and since $i$ preserve finite meets, $p, q \in \mathfrak{p}_i(m)
 \Rightarrow p \wedge q \in \mathfrak{p}_i(m)$.    
 
 Hence, for each $m \in \Sub{X}{\mathsf{M}}$, $\mathfrak{p}_i(m)$ is a filter on
 $X$, and \Arr{\mathfrak{p}_i}{\opp{\Sub{X}{\mathsf{M}}}}{\Fil{X}} is a
 preneighbourhood  on $X$. 
 
 Further:
 \[
 p \in \mathfrak{p}_i(p) \Leftrightarrow p \leq i(p) \Leftrightarrow i(p)
 = p \Rightarrow \mo_{\mathfrak{p}_i} = \bigl\{p \in \Sub{X}{\mathsf{M}}:
 i(p) = p\bigr\}
 \]
 shows for any $T \subseteq \mo_{\mathfrak{p}_i}$, $t \in T \Rightarrow t \leq
 \bigvee T \Rightarrow t = i(t) \leq i(\bigvee T) \Rightarrow \bigvee T 
 \leq i(\bigvee T)$, and hence $\mo_{\mathfrak{p}_i}$ is closed under arbitrary
 joins.   
 
 Moreover, for any $m \in \Sub{X}{\mathsf{M}}$:
 \begin{align*}
  \intr{\mathfrak{p}_i}{m} = \bigvee\bigl\{p \in \mo_{\mathfrak{p}_i}: p \leq
  m\bigr\} 
  = \bigvee\bigl\{p \in \Sub{X}{\mathsf{M}}: i(p) = p \leq
  m\bigr\} = i(m), \\
  p \in \mathfrak{p}_i(m) \Leftrightarrow p \leq i(m) \Leftrightarrow p
  \leq \intr{\mathfrak{p}_i}{m}, \\ \intertext{ and for any $S \subseteq
  \Sub{X}{\mathsf{M}}$}
  p \in \bigcap_{x \in S}\mathfrak{p}_i(x) \Leftrightarrow (\forall x \in
  S)\bigl(x \leq i(p)\bigr) \Leftrightarrow \bigvee S \leq i(p)
  \Leftrightarrow p \in \mathfrak{p}_i(\bigvee S).
 \end{align*}
 
 Hence $\mathfrak{p}_i \in \nhd{X} \subseteq 
 \low{\CAL{P}}{X}$. Furthermore, for any $\mu \in \CAL{P}_X$, if $i \leq
 \intr{\mu}{}$ then:
 \begin{equation*}
  p \in \mathfrak{p}_i(m) \Rightarrow m \leq i(p) \leq \intr{\mu}{p}
  \leq p \Rightarrow p \in \mu(m),
 \end{equation*}
 shows that the assignment $i \mapsto \mathfrak{p}_i$ extends to an
 order preserving map \Arr{\mathfrak{p}}{\CAL{K}_X}{\CAL{P}_X} such that
 \adjt{\mathfrak{p}}{\intr{\mu}{}} with $\intr{}{}\circ\mathfrak{p} =
 \id{\CAL{K}_X}$. 

  Clearly, from the adjunction the fibre \finv{\intr{}{}}{i} of any $i
 \in \CAL{K}_X$ has the neighbourhood $\mathfrak{p}_i$ as the smallest
 element. 
 
 Finally, if $\mu \in \CAL{P}_X$ be a neighbourhood on $X$ then from Theorem
 \ref{snbdgenbyopens}: 
 \begin{equation*}
  p \in \mu(m) \Leftrightarrow m \leq \intr{\mu}{p} \Leftrightarrow p
   \in \mathfrak{p}_{\intr{\mu}{}}(m), 
 \end{equation*}
 yields along with the observation $\mathfrak{p}_i$ for each $i \in
 \CAL{K}_X$ is a neighbourhood that $\mu \in \CAL{P}_X$ is a
 neighbourhood if and only if  $\mu = \mathfrak{p}_{\intr{\mu}{}}$,
 completing the proof. 
\end{proof} 

\subsubsection{Neighbourhoods and Pseudo-frame subsets}

\begin{Df}
 \label{pseudoframesets-df}
 A set $\mathfrak{O} \subseteq \Sub{X}{\mathsf{M}}$ of admissible
 subobjects of $X$ is said to be a {\em pseudo-frame set} if it is
 closed under finite meets and arbitrary joins.

 We denote the set of all pseudo-frame sets by \pfs{X} and is
 ordered by usual set inclusion.
\end{Df}

\begin{enumerate}[itemsep=1.2ex,label=\underline{Remark \arabic*},ref=\underline{Remark \arabic*},align=left,resume=rem]
 \item Clearly $\bigl\{0, \id{X}\bigr\}$ is the smallest and
       \Sub{X}{\mathsf{M}} is the largest element of \pfs{X}. \\

 \item Since any intersection of pseudo-frame sets is again a
       pseudo-frame set, it follows that \pfs{X} is a complete lattice
       with intersection being the arbitrary meet.

       Hence for every subset $\mathfrak{T} \subseteq \pfs{X}$ the
       supremum $\bigvee\mathfrak{T}$ exists, but a simple intrinsic
       description may be difficult to obtain. However, in case when
       \Sub{X}{\mathsf{M}} is itself a frame it has a simple description
       --- $\bigvee\mathfrak{T}$ is the set of all arbitrary joins of
       finite meets of elements of $\bigcup\mathfrak{T}$.
\end{enumerate}

\begin{Thm}
 \label{pfs=nhd}
 Given $\mathfrak{O} \in \pfs{X}$ let:
 \begin{equation}
  \label{pfs2nhd-eq}
   \mu_\mathfrak{O}(m) = \bigcup\bigl\{\atleast{q}: m \leq q \in
   \mathfrak{O}\bigr\}, m \in \Sub{X}{\mathsf{M}}.
 \end{equation}

 The assignment $\mathfrak{O} \mapsto \mu_\mathfrak{O}$ is an
 isomorphism of the complete lattices \pfs{X} and \nhd{X}.
\end{Thm}

\begin{proof}
 Since:
 \begin{itemize}
  \item $\mu_\mathfrak{O}(0) = \Sub{X}{\mathsf{M}}$,
	$\mu_\mathfrak{O}(\id{X}) = \bigl\{\id{X}\bigr\}$,
  \item $p \in \mu_\mathfrak{O}(m) \Rightarrow m \leq p$,
  \item $m \leq n \Rightarrow \bigl\{q \in \mathfrak{O}: q \geq n\bigr\}
  \subseteq \bigl\{q \in \mathfrak{O}: q \geq m\bigr\} \Rightarrow 
	\mu_\mathfrak{O}(m) \geq \mu_\mathfrak{O}(n)$,   
  \item $p' \geq p \in \mu_\mathfrak{O}(m) \Leftrightarrow (\exists q
	\in \mathfrak{O})(m \leq q \leq p \leq p') \Rightarrow p' \in
	\mu_\mathfrak{O}(m)$,
  \item $p, p' \in \mu_\mathfrak{O}(m) \Leftrightarrow (\exists q, q'
	\in \mathfrak{O})\bigl(m \leq q \leq p\text{ and }m \leq q' \leq
	p'\bigr) \Rightarrow m \leq q \wedge q' \leq p \wedge p'$, and
	since $\mathfrak{O}$ is closed under finite meets, $p \wedge p'
	\in \mu_\mathfrak{O}(m)$,  
 \end{itemize}
 it follows that
 \Arr{\mu_\mathfrak{O}}{\opp{\Sub{X}{\mathsf{M}}}}{\Fil{X}} is indeed a
 preneighbourhood on $X$.

 Further, $p \in \mu_\mathfrak{O}(p) \Leftrightarrow (\exists q \in
 \mathfrak{O})(p \leq q \leq p) \Leftrightarrow p \in \mathfrak{O}$, it
 follows that $\mo_{\mu_{\mathfrak{O}}} = \mathfrak{O}$, implying 
 $\mu_\mathfrak{O} \in \low{\CAL{P}}{X}$ by Theorem
 \ref{interiorisopen-alt} (page \pageref{interiorisopen-alt}). In
 particular, from Theorem \ref{Kuratowskiisretractofpnbdwithintopen}
 (page \pageref{Kuratowskiisretractofpnbdwithintopen}),
 \intr{\mu_\mathfrak{O}}{} is a Kuratowski interior operation.

 Using \eqref{pfs2nhd-eq} on the equivalence between
 \ref{opennbdsgenerate} \& \ref{nbdiffintcontain} in Theorem
 \ref{interiorisopen-alt} (see page \pageref{nbdiffintcontain})
 indicates 
 $\mu_\mathfrak{O} = \mathfrak{p}_{\intr{\mu_\mathfrak{O}}{}}$, and hence from
 Theorem \ref{Kuratowskiisretractofpnbdwithintopen} again,
 $\mu_\mathfrak{O} \in \nhd{X}$.

 If $\mathfrak{O}, \mathfrak{O'} \in \pfs{X}$ $\mathfrak{O} \subseteq
 \mathfrak{O'}$ then:
 \begin{equation*}
  p \in \mu_\mathfrak{O}(m) \Leftrightarrow (\exists q \in
  \mathfrak{O})(m \leq q \leq p) \Rightarrow (\exists q \in
  \mathfrak{O'})(m \leq q \leq p) \Leftrightarrow p \in \mu_\mathfrak{O'}(m)
 \end{equation*}
 implies $\mu_\mathfrak{O} \leq \mu_\mathfrak{O'}$; conversely, if
 $\mu_\mathfrak{O} \leq \mu_\mathfrak{O'}$ then using $\mathfrak{O} =
 \mo_{\mu_{\mathfrak{O}}}$ one obtains:
 \begin{equation*}
  p \in \mathfrak{O} \Leftrightarrow p \in \mu_\mathfrak{O}(p) \leq
   \mu_{\mathfrak{O'}}(p) \Rightarrow p \in \mathfrak{O'}.
 \end{equation*}

 Hence $\mathfrak{O} \subseteq \mathfrak{O'} \Leftrightarrow
 \mu_\mathfrak{O} \leq \mu_\mathfrak{O'}$ with $\mathfrak{O} = \mathfrak{O'}
 \Leftrightarrow \mu_{\mathfrak{O}} = \mu_\mathfrak{O'}$.

 Thus, the function \Arr{P}{\pfs{X}}{\nhd{X}} defined by
 $P(\mathfrak{O}) = \mu_\mathfrak{O}$ is an order preserving bijection
 with $P(\bigl\{0, \id{X}\bigr\}) = \nabla$ and $P(\Sub{X}{\mathsf{M}})
 = \atleast{}$.

 Now let $T \subseteq \nhd{X}$, $\mathfrak{O} = \sup\bigl\{\mo_\tau:
 \tau \in T\bigr\}$ and $\mathfrak{o} = \bigcap_{\tau \in T}\mo_\tau$;
 since \pfs{X} is a complete lattice $\mathfrak{o}, \mathfrak{O} \in
 \pfs{X}$. Then: 
 \begin{itemize}
  \item $\tau \in T \Rightarrow \mu_\mathfrak{o} \leq \tau \leq
	\mu_\mathfrak{O}$.
  \item If $\nu, \mu \in \nhd{X}$ such that $\tau \in T \Rightarrow \nu
	\leq \tau \leq \mu$ then $\tau \in T \Rightarrow \mo_\nu
	\subseteq \mo_\tau \subseteq \mo_\mu$. Hence using the
	definition of the supremum and infimum in \pfs{X}, $\mo_\nu
	\subseteq \mathfrak{o} \Leftrightarrow \nu \leq
	\mu_\mathfrak{o}$ and $\mathfrak{O} \subseteq \mo_\mu
	\Leftrightarrow \mu_\mathfrak{O} \leq \mu$.  
 \end{itemize} 

 Hence $\mu_\mathfrak{o} = \inf T$ and $\mu_\mathfrak{O} = \sup T$ in
 \nhd{X}, completing the proof. 
\end{proof}

\begin{enumerate}[itemsep=1.2ex,label=\underline{Remark \arabic*},ref=\underline{Remark \arabic*},align=left,resume=rem]
 \item \label{supinfofnbd}
       The proof also provides the route for computing the suprema or
       infima in \nhd{X}. Given a $T \subseteq \nhd{X}$ to obtain the
       suprema (respectively, infima) in \nhd{X}:\\
       \begin{enumerate}[itemsep=1.2ex,label=(\alph*),align=left]
	\item compute the suprema (respectively, infima) $\mu = \bigvee
	      T$ (respectively, $\mu = \bigwedge T =
	      \cannbd{\bigl(\bigcap 
	      T\bigr)}$) in \wnhd{X}, \\
	\item compute the set $\mo$ of $\mu$-open sets,\\
	\item the candidate for the suprema (respectively, infima) is
	      then $\mu_\mathfrak{O}$. 
       \end{enumerate} 

 \item Thus there are three ways to identify a neighbourhood on an
       object $X$, which is summarised in Table \ref{nbd-alt}. 
       
       \begin{table}
	\begin{adjustbox}{width=1\textwidth,center=\textwidth}
	 \begin{tabular}{||l||l|l|l||} \hline\hline
	  \diagbox[linecolor=blue,font=\itshape\itshape,linewidth=1pt]{From}{To} & \nhd{X} & \pfs{X} & \low{\CAL{K}}{X}
	  \\ \hline 
	  \nhd{X} ($\nu$) & $\nu$ & $\mo_\nu = \bigl\{p \in \Sub{X}{\mathsf{M}}: p
		  \in \nu(p) \bigr\}$ & $\intr{\nu}{m} = \bigvee\bigl\{p \in
		      \mo_\nu: p \leq m\bigr\}$\\ \hline
	  \pfs{X} ($\mathfrak{O}$) & $\mu_\mathfrak{O}(m) =
	      \bigcup\bigl\{\atleast{q}: m \leq q \in \mathfrak{O}\bigr\}$ &
		  $\mathfrak{O}$ & $\iota_\mathfrak{O} = \bigvee\bigl\{p \in
		      \mathfrak{O}: p \leq m\bigr\}$\\ \hline
	  \low{\CAL{K}}{X} ($i$) & $\mu_i(m) = \bigl\{p \in
	      \Sub{X}{\mathsf{M}}: m \leq 
	      i(p)\bigr\}$ & $\mathfrak{O}_i = \bigl\{p \in \Sub{X}{\mathsf{M}}:
		  i(p) = p\bigr\}$ & $i$ \\ \hline\hline
	  \multicolumn{3}{c}{} \\ 
	 \end{tabular}
	\end{adjustbox}
	
	\caption{Conversion between Facets of Neighbourhood}
	\label{nbd-alt}
       \end{table}
       
\end{enumerate} 

\subsubsection{Internal Topological Spaces}

\begin{Df}
 \label{internaltopology-df}
 A neighbourhood  $\mu$ is a {\em topology} on $X$ if $\mo_\mu$
 is a frame in the partial order of \Sub{X}{\mathsf{M}}.

 The set of all internal topologies on $X$ in denoted by \itop{X}.
\end{Df}

Clearly, $\nabla$ is an internal topology. However, \atleast{} is an
internal topology if and only if  \Sub{X}{\mathsf{M}} is itself a
frame.

Moreover,  the order isomorphism between \nhd{X} and
\pfs{X} indicate:
 \[
  \mu \leq \nu \in \itop{X} \Rightarrow \mu \in \itop{X}.
 \]

 Thus \itop{X} being a down set of \nhd{X} is a complete meet
 subsemilattice of \nhd{X}. The following is immediate:

\begin{Thm}
 \label{topiscomplete<>largesttop}
 \itop{X} is a complete sublattice of \nhd{X} if and only if  there
 exists a largest topology on $X$. 
\end{Thm}
\subsection{Morphisms of Neighbourhoods}

\begin{Df}
 \label{prenbdmorphisms-df}
 Given the internal preneighbourhood spaces \opair{X}{\mu} and \opair{Y}{\phi},
 a morphism \Arr{f}{X}{Y} of \Bb{A} is a {\em preneighbourhood morphism}
 if for every admissible subobject $n$ of $Y$:
 \begin{equation}
  p \in \phi(n) \Rightarrow \finv{f}{p} \in \mu(\finv{f}{n}).
 \end{equation}

 The symbol \Arr{f}{(X, \mu)}{(Y, \phi)} denotes 
 $f$ is a preneighbourhood morphism.
 
\end{Df}

Clearly the adjunctions in Theorem \ref{image-|preimage} (page
\pageref{image-|preimage}) and Theorem
\ref{changeofbase-for-upsets} (page \pageref{changeofbase-for-upsets})
easily suggest the following equivalent formulations of a
preneighbourhood morphism.

\begin{Thm}[{\cite[\S3]{HolgateSlapal2011}}]
 \label{prenbdmorphism-alt}
 Given the internal preneighbourhood spaces \opair{X}{\mu} and
 \opair{Y}{\phi} the following are equivalent for any morphism
 \Arr{f}{X}{Y} of \Bb{A}: 
 \begin{enumerate}[label=(\alph*),ref=(\alph*),align=left]
  \item $f$ is a preneighbourhood morphism.
  \item \label{pnbdmorphism-imginvfil}
	$\invfil{f}{\phi(n)} \subseteq \mu(\finv{f}{n})$, for every
	admissible subobject $n \in \Sub{Y}{\mathsf{M}}$.
  \item \label{pnbdmorphism-codsmaller}
	$\phi(n) \subseteq \imgfil{f}{\mu(\finv{f}{n})}$, for every
	admissible subobject $n \in \Sub{Y}{\mathsf{M}}$.
  \item \label{pnbdmorphism-domlarger}
	$\invfil{f}{\phi(\img{f}{m})} \subseteq \mu(m)$, for every
	admissible subobject $m \in \Sub{X}{\mathsf{M}}$.
 \end{enumerate}
\end{Thm}

\begin{enumerate}[itemsep=1.2ex,label=\underline{Remark \arabic*},ref=\underline{Remark \arabic*},align=left,resume=rem]
 \item Seen via diagrams, \Arr{f}{(X, \mu)}{(Y, \phi)} is a
       preneighbourhood morphism if and only if  the square below
       denotes a natural transformation from the composites of order
       preserving maps on the right to the composites of order
       preserving maps on the left, where the order preserving maps are
       considered as functors:  
       \[
       \xymatrixcolsep{1.2in}
       \xymatrixrowsep{1.2in}
       \xymatrix{
       {\opp{\Sub{X}{\mathsf{M}}}} \ar@{<-}@<1.2ex>[r]^{\finv{f}{}}
       \ar@<-1.2ex>[r]_{\img{f}{}} \ar@{}[r]|{\bot} \ar[d]_{\mu}
       \ar@2{}[dr]|{\rotatebox{-30}{$\Leftarrow$}} &
       {\opp{\Sub{Y}{\mathsf{M}}}} \ar[d]^{\phi} \\
       {\Fil{X}} \ar@{<-}@<1.2ex>[r]^{\invfil{f}{}}
       \ar@<-1.2ex>[r]_{\imgfil{f}{}} \ar@{}[r]|{\bot} & {\Fil{Y}}
       }.
       \]
\end{enumerate}

\subsubsection{Universal Weak Neighbourhoods}

\begin{Thm}
 \label{nbdreflinpnbd}
 Given a preneighbourhood morphism \Arr{f}{(X, \mu)}{(Y, \phi)} from
 the internal preneighbourhood space \opair{X}{\mu} to the internal weak
 neighbourhood space \opair{Y}{\phi}, the preneighbourhood morphism
 \Arr{f}{(X, \cannbd{\mu})}{(Y, 
 \phi)} is the unique morphism between internal weak
 neighbourhood spaces such that the diagram $\xymatrix{ {(X,
	\mu)} \ar[r]^{\id{X}} \ar[dr]_{f} & {(X, \cannbd{\mu})}   
 \ar@{.>}[d]^{ f} \\ & {(Y, \phi)} }$ commutes.  
\end{Thm} 

\begin{proof}  
 It is enough to show that
 \Arr{\invfil{f}{\phi(\img{f}{})}}{\opp{\Sub{X}{\mathsf{M}}}}{\Fil{X}}
 is a weak neighbourhood structure on $X$.

 Having shown the assertion above, since $f$ is already a
 preneighbourhood morphism, $\invfil{f}{\phi(\img{f}{})} \leq \mu$
 (using Theorem \ref{prenbdmorphism-alt}\ref{pnbdmorphism-domlarger}),
 and \cannbd{\mu} being the largest weak neighbourhood smaller than
 $\mu$ would then immediately yield $\invfil{f}{\phi(\img{f}{})} \leq
 \cannbd{\mu}$. Hence, \Arr{f}{(X, \cannbd{\mu})}{(Y, \phi)} would
 become a preneighbourhood morphism, completing the proof.  

 Towards the proof of the assertion: since $\phi$ is a
  weak neighbourhood structure, using Theorem \ref{pnbdisnbd<>interpolative}
  (page \pageref{pnbdisnbd<>interpolative}), \ref{filteredjoinoffilters}
 (page \pageref{filteredjoinoffilters}) and the
  adjunction \adjt{\invfil{f}{}}{\imgfil{f}{}} by Theorem
  \ref{changeofbase-for-upsets} (page \pageref{changeofbase-for-upsets}),
 for each admissible subobject $m \in \Sub{X}{\mathsf{M}}$:
 \begin{multline*}
  \invfil{f}{\phi(\img{f}{m})} =
  \invfil{f}{\bigl(\bigcup\bigl\{\phi(a): a \in
  \phi(\img{f}{m})\bigr\}\bigr)} \\
  =  \bigvee\bigl\{\invfil{f}{\phi(a)}: a \in \phi(\img{f}{m})\bigr\}
  = \bigcup\bigl\{\invfil{f}{\phi(a)}: a \in \phi(\img{f}{m})\bigr\}.  
 \end{multline*}

  Hence for each $p \in \invfil{f}{\phi(\img{f}{m})}$ there exists an
  $a \in \phi(\img{f}{m})$ such that $p \in \invfil{f}{\phi(a)}$.

  Since $\phi(a) \subseteq \phi(\img{f}{\finv{f}{a}}) \Rightarrow
  \invfil{f}{\phi(a)} \subseteq \invfil{f}{\phi(\img{f}{\finv{f}{a}})}$
  and $a \in \phi(\img{f}{m}) \Rightarrow \finv{f}{a} \in   
  \invfil{f}{\phi(\img{f}{m})}$, the statements $\finv{f}{a} \in
  \invfil{f}{\phi(\img{f}{m})}$ and $p \in
  \invfil{f}{\phi(\img{f}{\finv{f}{a}})}$ follow, showing 
  \invfil{f}{\phi(\img{f}{})} is interpolative, completing the proof
  using Theorem \ref{pnbdisnbd<>interpolative} (page
  \pageref{pnbdisnbd<>interpolative}).
\end{proof} 

\subsubsection{Universal Neighbourhoods}

Since the set \nhd{X} of all neighbourhood structures on $X$ is a
complete lattice (Theorem \ref{pfs=nhd}, page \pageref{pfs=nhd}), given
any weak neighbourhood structure $\mu$ on $X$ one has the largest
neighbourhood structure $\cansnbd{\mu} = \bigvee\bigl\{\nu \in \nhd{X}:
\nu \leq \mu\bigr\}$ on $X$ smaller than $\mu$. If \opair{X}{\mu} is an
internal weak neighbourhood space, \opair{Y}{\phi} is an internal
neighbourhood space, \Arr{f}{(X, \mu)}{(Y, \phi)} is a preneighbourhood
morphism such that $f$ has preimage preserve join property then for any
$S \subseteq \Sub{X}{\mathsf{M}}$:
\[
 \invfil{f}{\phi(\img{f}{\bigl(\bigvee S\bigr)})} =
 \invfil{f}{\phi\bigl(\bigvee_{s \in S}\img{f}{s}\bigr)} =
 \invfil{f}{\bigl(\bigcap_{s \in S}\phi(\img{f}s)\bigr)} = \bigcap_{s
 \in S}\invfil{f}{\phi(\img{f}s)}
\]
shows \invfil{f}{\phi(\img{f}{})} to be a neighbourhood structure on
$X$. Since $f$ is a preneighbourhood morphism,
$\invfil{f}{\phi(\img{f}{})} \leq \mu$ implies
$\invfil{f}{\phi(\img{f}{})} \leq \cansnbd{\mu}$, yielding:
\begin{Thm}
 \label{snbdreflinwnbd}
 Given a preneighbourhood morphism \Arr{f}{(X, \mu)}{(Y, \phi)} from
 the internal weak neighbourhood space \opair{X}{\mu} to the internal
 neighbourhood space \opair{Y}{\phi} where $f$ has the preimage preserve
 join property, the preneighbourhood morphism
 \Arr{f}{(X, \cansnbd{\mu})}{(Y, 
 \phi)} is the unique morphism between internal
 neighbourhood spaces such that the diagram:
 \[
 \xymatrix{
 {(X, \mu)} \ar[r]^{\id{X}} \ar[dr]_{f} & {(X, \cansnbd{\mu})}   
 \ar@{.>}[d]^{ f} \\
 & {(Y, \phi)} }
 \]
 commutes.  
\end{Thm}

\section{Categories of Neighbourhood Structures}
\label{NbdCategories}

\begin{Df}
 \label{internalspaces}
 The following categories are now stipulated.
 \begin{enumerate}[label=(\alph*),align=left,ref=(\alph*),itemsep=1.2ex]
  \item  \label{pretopologicalspaces}
	 \pre{\Bb{A}} is the category of internal preneighbourhood
	 spaces \opair{X}{\mu} and preneighbourhood morphisms
	 \Arr{f}{(X, \mu)}{(Y, \phi)}. 

  \item  \label{wnbdspaces}
	 \wNHD{\Bb{A}} is the full subcategory of \pre{\Bb{A}}
	 consisting of all internal weak neighbourhood spaces. 

  \item  \label{nbdspaces}
	 \NHD{\Bb{A}} is the subcategory of \wNHD{\Bb{A}} consisting of
	 internal neighbourhood spaces \opair{X}{\mu} and
	 preneighbourhood morphisms \Arr{f}{(X, \mu)}{(Y, \phi)} between
	 internal neighbourhood spaces where the morphism \Arr{f}{X}{Y}
	 of \Bb{A} has the preimage preserve joins  property.   

  \item  \label{topologicalspaces}
	 \Int{\Top}{\Bb{A}} is the full subcategory of \NHD{\Bb{A}}
	 consisting of internal topological spaces. 
 \end{enumerate}
\end{Df}

\subsection{Reflective Subcategories}

\begin{enumerate}[itemsep=1.2ex,label=\underline{Remark \arabic*},ref=\underline{Remark \arabic*},align=left,resume=rem]
 \item \label{wnbdbireflectiveinpretop}
       Theorem \ref{nbdreflinpnbd} (page \pageref{nbdreflinpnbd})
       exactly shows \wNHD{\Bb{A}} to be a bireflective full subcategory
       of \pre{\Bb{A}}.

 \item \label{nbdbireflectiveinwnbd-ppj}
        Let \PPJ{\wNHD{\Bb{A}}} be the non-full subcategory of
       \wNHD{\Bb{A}} with internal weak neighbourhood spaces as objects
       and preneighbourhood morphisms \Arr{f}{(X, \mu)}{(Y, \phi)} for
       which the morphism \Arr{f}{X}{Y} of \Bb{A} have the preimage
       preserve join property. \\
 
       Theorem \ref{snbdreflinwnbd} (page \pageref{snbdreflinwnbd})
       shows the category \NHD{\Bb{A}} is a bireflective full
       subcategory of \PPJ{\wNHD{\Bb{A}}}.       
\end{enumerate}

\begin{Thm}
 \label{largesttop-alt}
 Let \PPJ{\Bb{A}} be the subcategory of \Bb{A} consisting of all objects
 of \Bb{A} and morphisms \Arr{f}{X}{Y} of \Bb{A} which have the preimage
 preserve joins property.
 
 Then, the following are equivalent:
 \begin{enumerate}[label=(\alph*),ref=(\alph*),align=left]
  \item \label{largest} For every object $X$, there exists a largest
	internal topological structure on $X$.
  \item \label{topreflinstrnbd} \Int{\Top}{\Bb{A}} is a full bireflective
	subcategory of \NHD{\Bb{A}}.
  \item \label{toptopological} \Int{\Top}{\Bb{A}} is topological over
	\InvRAdjt{\Bb{A}}. 
 \end{enumerate}
\end{Thm}

\begin{proof}
  
 \begin{description}
  \item[\underline{\ref{largest} implies \ref{topreflinstrnbd}}] Choose
	     and fix an internal neighbourhood space \opair{X}{\mu}.

	     Since from the assumption \itop{X} is a complete sublattice
	     of \nhd{X} (Theorem \ref{topiscomplete<>largesttop}, page
	     \pageref{topiscomplete<>largesttop}), $\cantop{\mu} =
	     \bigvee\bigl\{\nu \in \itop{X}: \nu \leq \mu\bigr\}$ is the
	     largest internal topology on $X$ smaller than $\mu$. Hence,
	     \Arr{\id{X}}{(X, \mu)}{(X, \cantop{\mu})} is a bimorphism
	     of \NHD{\Bb{A}}. 

	     If \Arr{f}{(X, \mu)}{(Y, \phi)} is a morphism of
	     \NHD{\Bb{A}} from the internal neighbourhood space
	     \opair{X}{\mu} to the internal topological space
	     \opair{Y}{\phi} then arguments similar to the one just
	     before the statement of Theorem \ref{snbdreflinwnbd} (page
	     \pageref{snbdreflinwnbd}) show \invfil{f}{\phi(\img{f}{})}
	     to be a neighbourhood structure on $X$. Furthermore using
	     \adjt{\img{f}{}}{\adjt{\finv{f}{}}{\fInv{f}{}}} (Theorem
	     \ref{image-|preimage} (page \pageref{image-|preimage}) \&
	     Theorem \ref{ppj full or
	     partial}\ref{preimagehasrightadjoint-alt} (page 
	     \pageref{preimagehasrightadjoint-alt})): 
	     \begin{multline*}
	      p \in \invfil{f}{\phi(\img{f}{p})}
	      \Leftrightarrow (\exists q \in
	      \phi(\img{f}{p}))(\finv{f}{q} \leq p)  
	      \Leftrightarrow (\exists q \in \phi(\img{f}{p}))(q \leq
	      \fInv{f}{p}) 
	      \Leftrightarrow \fInv{f}{p} \in \phi(\img{f}{p}) \\
	      \Leftrightarrow (\exists u \in \mo_\phi)(\img{f}{p} \leq u
	      \leq \fInv{f}{p})
	      \Leftrightarrow (\exists u \in \mo_\phi)(p \leq
	      \finv{f}{u} \leq p)
	      \Leftrightarrow (\exists u \in \mo_\phi)(p =
	      \finv{f}{u}),
	     \end{multline*}
	     shows:
	     \[
	     \mo_{\invfil{f}{\phi(\img{f}{})}} = \bigl\{\finv{f}{u}: u
	     \in \mo_\phi\bigr\}.  
	     \]

	     Since $\mo_\phi$ is a frame and \finv{f}{} preserves all
	     joins and meets,
	     $\mo_{\invfil{f}{\phi(\img{f}{})}}$ is also a frame. Hence
	     \invfil{f}{\phi(\img{f}{})} is an internal topology on $X$,
	     smaller than $\mu$, implying $\invfil{f}{\phi(\img{f}{})}
	     \leq \cantop{\mu}$. 

	     Consequently, \Arr{f}{(X, \cantop{\mu})}{(Y, \phi)} is the
	     unique morphism of internal topological spaces such that
	     the diagram $\xymatrix{ {(X, \mu)} \ar[r]^{\id{X}}
	     \ar[dr]_f & {(X, \cantop{\mu})} \ar@{.>}[d]^{f} \\ & {(Y,
	     \phi)} }$ commutes in \NHD{\Bb{A}}.  
	     
  \item[\underline{\ref{largest} implies \ref{toptopological}}]
	     Choose and fix a family $\bigl((X_i, \mu_i)\bigr)_{i \in 
	     I}$ of internal topological spaces and a family
	     $\bigl(\Arr{f_i}{X}{X_i}\bigr)_{i \in I}$ of morphisms 
	     from \PPJ{\Bb{A}}.

	     Since for each $i \in I$, $\invfil{f_i}{\mu_i(\img{f_i}{})}
	     \in \itop{X}$ and from our assumption \itop{X} is a
	     complete lattice, $\mu = \bigvee_{i \in
	     I}\invfil{f_i}{\mu_i(\img{f_i}{})} \in \itop{X}$. Hence
	     \Arr{f_i}{(X, \mu)}{(X_i, \mu_i)}, for each $i \in I$, is a
	     morphism of \Int{\Top}{\Bb{A}}.

	     If \opair{Z}{\zeta} be an internal topological space and
	     \Arr{g}{Z}{X} be a morphism of \PPJ{\Bb{A}} such that for
	     each $i \in I$, \Arr{f_i\circ g}{(Z, \zeta)}{(X_i, \mu_i)}
	     is a morphism of \Int{\Top}{\Bb{A}}, then for each $z \in
	     \Sub{Z}{\mathsf{M}}$, $\zeta(z) \supseteq
	     \overleftarrow{(f_i\circ g)}\mu_i(\img{(f_i\circ g)}{z}) =
	     \invfil{g}{\bigl(\invfil{f_i}{\mu_i(\img{f_i}{(\img{g}{z})})}\bigr)}
	     \Leftrightarrow \imgfil{g}{\zeta(z)} \supseteq
	     \invfil{f_i}{\mu_i\bigl(\img{f_i}{(\img{g}{z})}\bigr)}$.

	     Hence, for each $z \in \Sub{Z}{\mathsf{M}}$,
	     $\mu(\img{g}{z})  
	     \subseteq \imgfil{g}{\zeta(z)} \Leftrightarrow 
	     \invfil{g}{\mu(\img{g}{z})} \subseteq \zeta(z)$, implying
	     \Arr{g}{(Z, \zeta)}{(X, \mu)} to be a morphism of internal
	     topological spaces, and the unique one making each
	     $\xymatrix{ {(X, \mu)} \ar@{<.}[dr]_{!\, g} \ar[r]^{f_i} &
	     {(X_i, \mu_i)} \ar@{<-}[d]^{f_i\circ g} \\ & {(Z, \zeta)}
	     }$ ($i \in I$) to commute, proving \ref{toptopological}.  

  \item[\underline{\ref{topreflinstrnbd} implies \ref{largest}}]
	      Assuming \ref{topreflinstrnbd}, given any internal
	     topology $\mu$ on $X$, one has the diagram $\xymatrix{ {(X,
	     \atleast{})} \ar[r]^{\id{X}} \ar[dr]_{\id{X}} & {(X,
	     \overline{\atleast{\,}})} \ar@{.>}[d]^{!\,\id{X}} \\ & {(X,
	     \mu)} }$ to commute uniquely, yielding $\mu \leq 
	     \overline{\atleast{}}$, proving \ref{largest}.

  \item[\underline{\ref{toptopological} implies \ref{largest}}]
	      Assuming \ref{toptopological}, consider the object $X$,
	     the family $\bigl((X, \mu)\bigr)_{\mu \in \tnhd{X}}$ of
	     internal topological objects and the family
	     $\bigl(\id{X}\bigr)_{\mu \in \tnhd{X}}$ of morphisms of
	     \InvRAdjt{\Bb{A}}. From hypothesis, there exists a unique
	     internal topological structure $\lambda_X$ on $X$ such that
	     for each $\mu \in \tnhd{X}$, \Arr{\id{X}}{(X,
	     \lambda_X)}{(X, \mu)} is a morphism of \itop{X}, implying
	     $\mu \leq \lambda_X$, for each $\mu \in \tnhd{X}$, proving
	     \ref{largest}.   
 \end{description}

\end{proof}

\subsection{Results on Topologicity}

\begin{Thm}
 \label{topologicityresults}
 \begin{enumerate}[label=(\alph*),ref=(\alph*),align=left]
  \item \label{pnbdtopoverbase} The category \pNHD{\Bb{A}} of internal
	preneighbourhood spaces is topological over \Bb{A}.
  \item \label{wnbdtopoverbase} The category \wNHD{\Bb{A}} of internal
	preneighbourhood spaces is topological over \Bb{A}. 
  \item \label{nbdtopoverppj} The category \NHD{\Bb{A}} of internal
	neighbourhood spaces is topological over \PPJ{\Bb{A}}. 	
 \end{enumerate}
\end{Thm}

\begin{proof}
 The proof follows from the facts:
 \begin{itemize}
  \item if \Arr{f}{X}{Y} is a morphism of \Bb{A} and \opair{Y}{\phi} is
	an internal preneighbourhood space (respectively, an internal
	weak neighbourhood space), then from the proof of Theorem
	\ref{nbdreflinpnbd} (page \pageref{nbdreflinpnbd})
	\invfil{f}{\phi(\img{f}{})} is a preneighbourhood (respectively,
	weak neighbourhood) structure on $X$, 
  \item if \Arr{f}{X}{Y} is a morphism from \PPJ{\Bb{A}} and
	\opair{Y}{\phi} is an internal neighbourhood space then using
	arguments just before Theorem \ref{snbdreflinwnbd} (page
	\pageref{snbdreflinwnbd}) \invfil{f}{\phi(\img{f}{})} is a
	neighbourhood structure on $X$, and
  \item the sets \pnhd{X}, \wnhd{X} and \nhd{X} of preneighbourhood
	structures, weak neighbourhood structures and neighbourhood
	structures on $X$ make a complete lattice --- Theorem
	\ref{nbdstructuresmakecompletelattice} (page
	\pageref{nbdstructuresmakecompletelattice}), Theorem
	\ref{pfs=nhd} (page \pageref{pfs=nhd}) and 
	\ref{supinfofnbd} (page \pageref{supinfofnbd}). 
 \end{itemize}
\end{proof}

All of this leads to the diagram in Figure \ref{catstate} (page
\pageref{catstate}), which summarises the results obtained so far. While
the general situation appears in Figure \ref{catstate-1}, the picture is
simplified when every morphism has preimage preserve join property (see
Figure \ref{catstate-2}). As a consequence of Corollary
\ref{ppj>subobjectframes} (page \pageref{ppj>subobjectframes}) every
lattice of admissible subobjects is a frame. This is the situation for
$\Bb{A} = \Set$, in particular.

\begin{center}
 \begin{figure}[p]
  \begin{subfigure}[p]{0.96\linewidth}
   \xymatrixcolsep{1.2in}
   \xymatrixrowsep{1.2in}
   \xymatrix{
   {\pre{\Bb{A}}} \ar[ddr]|(0.36){U}|(0.48){\text{\tiny Theorem
   \ref{topologicityresults}\ref{pnbdtopoverbase}\hspace{2ex}}}
   \ar@{<-}[r]^{\text{\tiny bireflective subcategory}}_{\text{\tiny
   Theorem \ref{nbdreflinpnbd}}} & {\wNHD{\Bb{A}}}
   \ar@{<-}[r]^{\text{\tiny non full subcategory}}
   \ar[dd]|(0.36){V}|(0.48){\text{\tiny Theorem
   \ref{topologicityresults}\ref{wnbdtopoverbase}}} &
   {\PPJ{\wNHD{\Bb{A}}}} \ar@{<-}[r]^{\text{\tiny bireflective
   subcategory}}_{\text{\tiny Theorem \ref{snbdreflinwnbd}}} \ar[dd] &
   {\NHD{\Bb{A}}} \ar@{<.}[d]|(0.48){\text{\tiny bireflective
   subcategory}}|(0.66){\text{\tiny Theorem \ref{largesttop-alt}}}
   \ar@/_3ex/[ddl]|(0.48){W}|(0.6){\text{\tiny Theorem
   \ref{topologicityresults}\ref{nbdtopoverppj}}} &  \\ 
   & & & {\Int{\Top}{\Bb{A}}} \ar@{.>}[dl]^{T}|(0.36){\text{\tiny
   Theorem \ref{largesttop-alt}\hspace{3ex}}} \\ 
   & {\Bb{A}} \ar@{<-}[r]_{\text{non full subcategory}} & {\PPJ{\Bb{A}}}
   }
   
   \caption{Categories of Neighbourhood Structures: the {\em
   forgetful} functors $U$, $V$, $W$ and $T$ (modulo its existence)
   are topological.  The dotted lines are used to highlight extra
   necessary conditions.}
   \label{catstate-1}   
  \end{subfigure}
  \vskip\baselineskip
  
  \begin{subfigure}[p]{\linewidth}
   \xymatrixcolsep{1.44in}
   \xymatrixrowsep{1.2in}
   \xymatrix{
   {\pre{\Bb{A}}} \ar@/_3ex/[dr]|(0.36){U}|(0.48){\text{Theorem
   \ref{topologicityresults}\ref{pnbdtopoverbase}}}
   \ar@{<-}[r]^{\text{\tiny bireflective subcategory}}_{\text{Theorem
   \ref{nbdreflinpnbd}}} & {\wNHD{\Bb{A}}} \ar@{<-}[r]^{\text{\tiny
   bireflective subcategory}}_{\text{Theorem \ref{snbdreflinwnbd}}}
   \ar[d]|(0.36){V}|(0.54){\text{Theorem
   \ref{topologicityresults}\ref{wnbdtopoverbase}}} 
   & {\NHD{\Bb{A}}} \ar[dl]|(0.36){W}|(0.54){\text{Theorem
   \ref{topologicityresults}\ref{nbdtopoverppj}}}
   \ar@2{-}[d] \\
   & {\Bb{A}} \ar@{<-}[r]_T & {\Int{\Top}{\Bb{A}}}
   }

   \caption{Categories of Neighbourhood Structures when every morphism
   has preimage preserve join property, as a consequence of which every
   lattice of admissible subobjects is a frame.}
   \label{catstate-2}   
  \end{subfigure}
  
  \vskip\baselineskip
  \caption{Summarising Categories of Neighbourhood Structures}
  \label{catstate}  
 \end{figure}
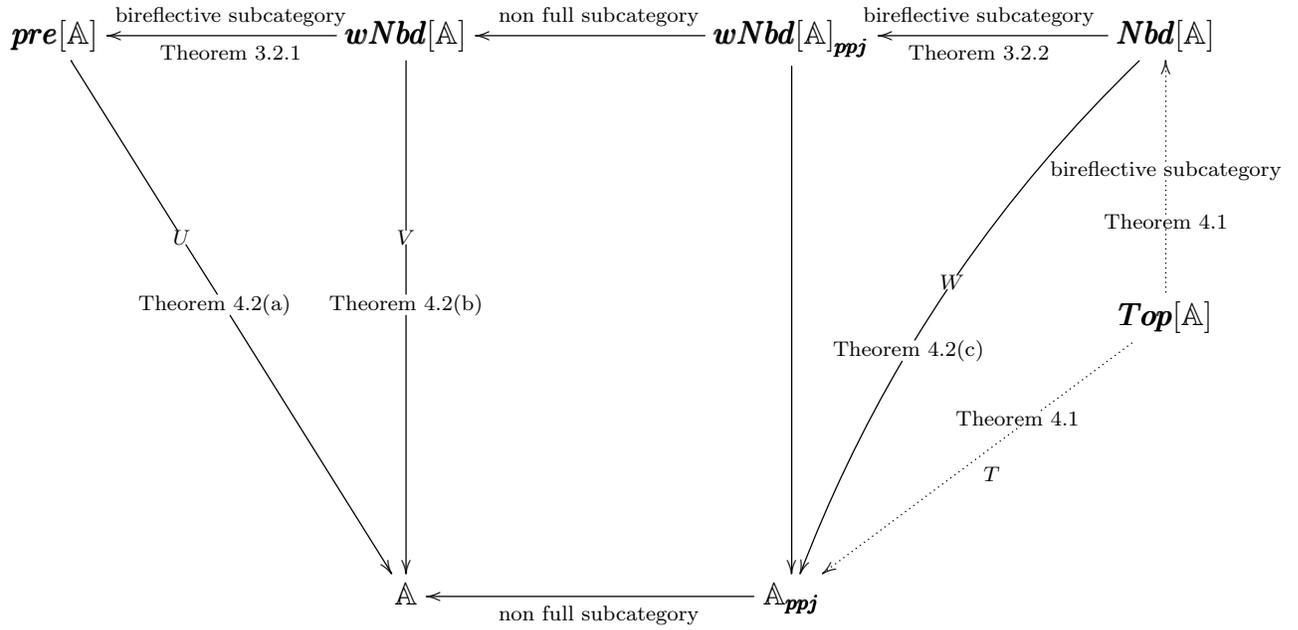
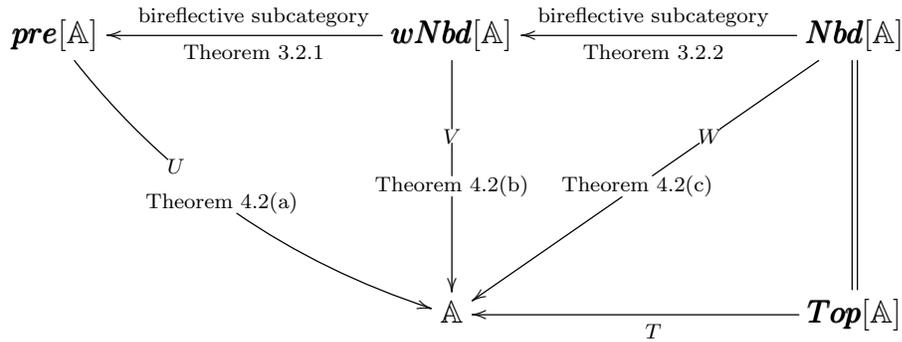
\end{center}

\section{Regular Epimorphisms of Internal Neighbourhood Spaces}
\label{RegEpi}

\subsection{Regular Epimorphisms of \pre{\Bb{A}}}

\begin{Thm}
 \label{regepipreobj-desc}
 A morphism \Arr{f}{(X, {\gamma})}{(Y, {\phi})} of \pre{\Bb{A}} is a  
 regular epimorphism if and only if  the morphism
 \Arr{f}{X}{Y} is a regular epimorphism of \Bb{A} and: 
 \begin{equation}
  \label{regepipretop-eq}
   \phi(y) = \bigl\{u \in \Sub{Y}{\mathsf{M}}: y \leq u\text{ and
   }\finv{f}u \in \gamma(\finv{f}y)\bigr\}.
 \end{equation}
\end{Thm}

\begin{proof}
 \hfill{\,}
  
 \begin{itemize}
  \item[\underline{if} part:\hspace{5.5ex}] Since \Arr{f}{X}{Y} is a
	     regular epimorphism of \Bb{A}, it is the coequaliser of its
	     kernel pair $\xymatrix{ {\Kerp{f}} \ar[r]^(0.6){p_2}
	     \ar[d]_{p_1} & {X} \ar[d]^f \\ {X} \ar[r]_f & {Y} }$.

	     Since the forgetful functor \Arr{U}{\pre{\Bb{A}}}{\Bb{A}}
	     creates kernel pairs, there exists a unique
	     pre-neighbourhood $\kappa$ on \Kerp{f} such that
	     $\xymatrix{ {(\Kerp{f}, \kappa)} \ar[r]^(0.6){p_2}
	     \ar[d]_{p_1} & {(X, \gamma)} \ar[d]^f \\ {(X,
	     \gamma)} \ar[r]_f & {(Y, \phi)} }$ is the
	     kernel pair of \Arr{f}{(X, \gamma)}{(Y,
	     \phi)} in \pre{\Bb{A}}.

	     Let \Arr{g}{(X, {\gamma})}{(Z, \zeta)} be a
	     pretopological morphism such that $\comp{g}{p_1} =
	     \comp{g}p_2$. Then:
	     \begin{itemize}
	      \item From the coequaliser in \Bb{A}:
		    \[
		    \xymatrix{
		    {\Kerp{f}} \ar@<0.6ex>[r]^(0.6){p_1}
		    \ar@<-0.6ex>[r]_(0.6){p_2} & {X} 
		    \ar[r]^f \ar[dr]_g & {Y} \ar@{.>}[d]^{!\, h} \\
		    & & {Z}
		    },
		    \]
		    there exists the unique morphism \Arr{h}{Y}{Z} such
		    that $g = \comp{h}{f}$. 
	      \item Choose and fix admissible subobjects $u, z$ of $Z$
		    with $u \in \zeta(z)$.

		    Since \Arr{g}{(X, \gamma)}{(Z, \zeta)} is a
		    pre-neighbourhood morphism, $\finv{g}{u} \in
		    \gamma(\finv{g}z)$.

		    But $u \in \zeta(z) \Rightarrow z \leq u \Rightarrow
		    \finv{h}z \leq \finv{h}u$ and:
		    \begin{equation*}
		     \finv{g}{u} \in \gamma(\finv{g}z) \Leftrightarrow
		     \finv{(h\circ f)}u \in \gamma(\finv{(h\circ f)}{z})
		     \\
		     \Leftrightarrow \finv{f}{(\finv{h}u)} \in
		     \gamma(\finv{f}{(\finv{h}z)}),
		    \end{equation*}
		    so that:
		    \[
		     u \in \zeta(z) \Rightarrow \finv{h}z \leq
		    \finv{h}u\text{ and }\finv{f}{(\finv{h}u)} \in
		    \gamma(\finv{f}{(\finv{h}z)}) \Leftrightarrow
		    \finv{h}u \in \phi(\finv{h}z),
		    \]
		    proving \Arr{h}{(Y, \phi)}{(Z, \zeta)} to be a
		    pre-neighbourhood morphism.
	      \item Since $U$ is faithful, \Arr{h}{(Y,
		    \phi)}{(Z, \zeta)} is the unique
		    pretopological morphism which makes the diagram:
		    \[
		    \xymatrix{
		    {(\Kerp{f}, \kappa)}
		    \ar@<0.6ex>[r]^(0.6){p_1} \ar@<-0.6ex>[r]_(0.6){p_2}
		    & {(X, \gamma)} \ar[r]^f \ar[dr]_g & {(Y,
		    \phi)} \ar@{.>}[d]^{!\, h} \\  
		    & & {(Z, \zeta)} 
		    }
		    \]
		    to commute in \pre{\Bb{A}}.
	     \end{itemize}

	     Hence \Arr{f}{(X, \gamma)}{(Y, \phi)} is a
	     regular epimorphism in \pre{\Bb{A}}.
  \item[\underline{only if} part:] Since the forgetful functor
	     \Arr{U}{\pre{\Bb{A}}}{\Bb{A}} preserve coequalisers, the
	     coqualiser diagram:
	     \[
	     \xymatrix{
	     (Z, \zeta) \ar@<0.6ex>[r]^p \ar@<-0.6ex>[r]_q &
	     {(X, \gamma)} \ar[r]^f & {(Y, \phi)} 
	     }
	     \]
	     in \pre{\Bb{A}} is mapped to the coequaliser diagram:
	     \[
	     \xymatrix{
	     {Z} \ar@<0.6ex>[r]^p \ar@<-0.6ex>[r]_q & {X}
	     \ar[r]^f & {Y}
	     }
	     \]
	     in \Bb{A}. In particular, $f$ is a regular epimorphism of
	     \Bb{A}, and since \fact{\mathsf{E}}{\mathsf{M}} is
	     proper, $\RegEpi{\Bb{A}} \subseteq \ExtEpi{\Bb{A}}
	     \subseteq \mathsf{E}$, implying $f \in \mathsf{E}$.

	     Define:
	     \begin{equation}
	      \label{selectingsubinvinpretop-eq}
	      \psi(y) = \biggl\{v \in \Sub{Y}{\mathsf{M}}: y \leq
	      v\text{ and }\finv{f}v \in \gamma(\finv{f}y)\biggr\},
	      \quad\text{ for all }y \in \Sub{Y}{\mathsf{M}}.  
	     \end{equation}

	     Then:
	     \begin{itemize}
	      \item $u \geq v \in \psi(y)$ implies $u \geq v \geq y$ and 
		    $\finv{f}u \geq \finv{f}v \in \gamma(\finv{f}y)
		    \Rightarrow \finv{f}u \in \gamma(\finv{f}y)$, 
		    since $\gamma(\finv{f}y)$ is a filter.

		    Hence $u \in \psi(y)$, showing $\psi(y)$ is an 
		    upset. 
	      \item $u, v \in \psi(y)$ implies $u, v \geq y$ and
		    $\finv{f}u, \finv{f}v \in \gamma(\finv{f}y)$,
		    implying $u \wedge v \geq y$ and $\finv{f}{(u \wedge
		    v)} = \finv{f}u \wedge \finv{f}v \in
		    \gamma(\finv{f}y)$. Hence $u \wedge v \in \psi(y)$,
		    showing $\psi(y)$ to be closed under finite 
		    intersections. 
	      \item $y \geq z$ and $v \in \psi(y)$ implies $v \geq y
		    \geq z$ and $\finv{f}v \in \gamma(\finv{f}y)
		    \subseteq \gamma(\finv{f}z)$, so that in particular,
		    $v \geq z$ and $\finv{f}v \in \gamma(\finv{f}z)$,
		    i.e., $v \in \psi(z)$. Hence $\psi(y) \subseteq  
		    \psi(z)$, showing the assignment $y \mapsto \psi(y)$
		    to be an order preserving map
		    \Arr{\psi}{\opp{\Sub{Y}{\mathsf{M}}}}{\Fil{Y}}.
	      \item Since $v \in \psi(y)$ implies $v \geq y$, it
		    follows that $\psi$ is a pre-neighbourhood on $Y$.
	      \item Since $f$ is a pretopological morphism, if $u \in
		    \phi(y)$ then $\finv{f}{u} \in
		    \gamma(\finv{f}{y})$. Hence, $u \geq y$ and
		    $\finv{f}u \in \gamma(\finv{f}y)$ implies $u \in
		    \psi(y)$. Hence $\phi(y) \subseteq \psi(y)$, as 
		    filters and hence $\phi \leq \psi$, as
		    pre-neighbourhoods on $Y$. 
	      \item From the very definition of $\psi$, for any $y \in
		    \Sub{Y}{\mathsf{M}}$ and a $v \in \psi(y)$, 
		    \finv{f}{v} is in $\gamma(\finv{f}{y})$, so
		    that \Arr{f}{(X, \gamma)}{(Y, \psi)} is also a
		    pretopological morphism.
	      \item Since \Arr{f}{(X, \gamma)}{(Y,
		    \phi)} is a regular epimorphism in
		    \pre{\Bb{A}}, 
		    there exists a unique pretopological morphism
		    \Arr{h}{(Y, \phi)}{(Y, \psi)} making the diagram:
		    \[
		    \xymatrix{
		    (Z, \zeta) \ar@<0.6ex>[r]^p \ar@<-0.6ex>[r]_q & {(X,
		    \gamma)} \ar[r]^f \ar[dr]_f & {(Y, \phi)} 
		    \ar@{.>}[d]^{!\, h} \\
		    & & {(Y, \psi)}}
		    \]
		    to commute in \pre{\Bb{A}}.		    	
	      \item Since $f \in \RegEpi{\Bb{E}}$, 
		    $h\circ f = f = \id{Y}\circ f$ in \Bb{A} implies $h
		    =  \id{Y}$. 
	      \item Hence \Arr{\id{Y}}{(Y, \phi)}{(Y, \psi)} is a
		    pre-neighbourhood morphism, entailing $\psi \leq
		    \phi$. 
	     \end{itemize}

	     Hence $\phi = \psi$, completing the proof.
 \end{itemize}
\end{proof}

\begin{enumerate}[label=\underline{Remark \arabic*},ref=\underline{Remark \arabic*},align=left,resume=rem,itemsep=1.2ex]
 \item The proof only requires the forgetful functor
       \Arr{U}{\pre{\Bb{A}}}{\Bb{A}} to create kernel pairs and preserve
       coequalisers. Theorem
       \ref{topologicityresults}\ref{pnbdtopoverbase} provides much more
       than just this requirement.
 \item \label{extracondition}
       If $\finv{f}{u} \in \gamma(\finv{f}{y})$ then $\finv{f}{y} \leq
       \finv{f}{u} \Leftrightarrow \img{f}{\finv{f}{y}} \leq u$, and
       hence the extra restriction in \eqref{regepipretop-eq} on page
       \pageref{regepipretop-eq} is to ensure
       the description of regular epimorphisms in \pre{\Bb{A}}.
 \item \label{pullbackstability>simpledescregepi}
       The condition of a regular epimorphism $f$ being stably in
       $\mathsf{E}$ is necessary to ensure simpler description of the
       regular epimorphism \Arr{f}{(X, \gamma)}{(Y, \phi)}:
       \[
	u \in \phi(y) \Leftrightarrow \finv{f}{u} \in \gamma(\finv{f}{y}).
       \]

       This is exactly the situation in case when $\Bb{A} = \Set$. 
\end{enumerate}

\subsection{Regular Epimorphisms of \NHD{\Bb{A}}}

Given a morphism \Arr{f}{(X, \gamma)}{(Y, \phi)} of preneighbourhoods,
the proof of Theorem \ref{regepipreobj-desc} suggests:
\[
 \phi(y) = \bigl\{u \in \Sub{Y}{\mathsf{M}}: y \leq u\text{ and
 }\finv{f}u \in \gamma(\finv{f}y)\bigr\},
\]
is actually a preneighbourhood on $Y$, $\phi \leq \psi$ and \Arr{f}{(X,
\gamma)}{(Y, \psi)} is a morphism of preneighbourhoods.

Further:
\[
 \mo_\psi = \bigl\{u \in \Sub{Y}{\mathsf{M}}: \finv{f}{u} \in
 \mo_\gamma\bigr\},   
\]
is closed under arbitrary joins, if $\gamma$ preserve arbitrary meets.

Hence, if $\gamma$ is preneighbourhood which preserve arbitrary meets
then \intr{\psi}{} is a Kuratowski operator. Consequently the smallest
preneighbourhood $\hat\psi$ in the fibre
\finv{\intr{}{}}{\bigl(\intr{\psi}{}\bigr)} of \intr{\psi}{} is a
neighbourhood on $Y$ (see Theorem
\ref{Kuratowskiisretractofpnbdwithintopen}, page
\pageref{Kuratowskiisretractofpnbdwithintopen}). Hence $\mo_{\hat\psi} =
\mo_{\psi} \supseteq \mo_\phi$ and $\hat\psi \leq \psi$. If, further
$\phi$ be a neighbourhood then using Theorem \ref{pfs=nhd} (see page
\pageref{pfs=nhd}) $\phi \leq \hat\psi \leq \psi$. All these
observations along with the topologicity of \NHD{\Bb{A}} over
\PPJ{\Bb{A}} (\ref{nbdbireflectiveinwnbd-ppj}, page
\pageref{nbdbireflectiveinwnbd-ppj}) yield similarly as in
Theorem \ref{regepipreobj-desc}:

\begin{Thm}
 \label{regepinbdobj-desc}
 A morphism \Arr{f}{(X, {\gamma})}{(Y, {\phi})} of \NHD{\Bb{A}} is a  
 regular epimorphism if and only if  the morphism
 \Arr{f}{X}{Y} is a regular epimorphism of \PPJ{\Bb{A}} and: 
 \begin{equation}
  \label{regepinbd-eq}
   \phi(y) = \bigl\{u \in \Sub{Y}{\mathsf{M}}: y \leq u\text{ and
   }\finv{f}u \in \gamma(\finv{f}y)\bigr\}.
 \end{equation}
\end{Thm}

\begin{enumerate}[label=\underline{Remark \arabic*},ref=\underline{Remark \arabic*},align=left,resume=rem,itemsep=1.2ex]
 \item In case where $\gamma$ is a neighbourhood and \finv{f}{} preserve
       arbitrary joins then for any $S \subseteq \Sub{Y}{M}$:
\begin{multline*}
 u \in \psi\bigl(\bigvee S\bigr) \Leftrightarrow u \geq \bigvee S\text{
 and }\finv{f}{u} \in \gamma\biggl(\finv{f}{\bigl(\bigvee
 S\bigr)}\biggr) = \gamma\bigl(\bigvee_{s \in S}\finv{f}{s}\bigr)\\
 \Leftrightarrow \bigvee S \leq u \text{ and }\finv{f}{u} \in
 \bigcap_{s \in S}\gamma(\finv{f}{s}) \\
 \Leftrightarrow (\forall s \in S)\bigl(s \leq u \text{ and }
 \finv{f}{u} \in \gamma(s)\bigr) \\
 \Leftrightarrow (\forall s \in S)\bigl(u \in \psi(s)\bigr)
 \Leftrightarrow u \in \bigcap_{s \in S}\psi(s),
\end{multline*}
shows $\psi$ to preserve meets. However, this does not guarantee
       whether $\psi$ is a neighbourhood.
\end{enumerate}

\subsection{Hereditary Regular Epimorphisms}

Given any preneighbourhood $\gamma$ of an object $X$ of \Bb{A} and an
admissible subobject \Arr{p}{P}{X} of $X$ there exists from
topologicity of the forgetful functor \Arr{U}{\pre{\Bb{A}}}{\Bb{A}} a
unique {\em smallest} preneighbourhood $\gamma_p$ on $P$ such that
\Arr{p}{(P, \gamma_p)}{(X, \gamma)} is a preneighbourhood
morphism. Indeed:
\begin{equation}
 \label{subobjectprenbd-eq}
  \gamma_p(m) = \bigl\{u \in \Sub{P}{\mathsf{M}}: (\exists w \in
  \gamma(p\circ m))\bigl(p \wedge w \leq p\circ u\bigr)\bigr\}, \quad\text{ for
  all }m \in \Sub{U}{\mathsf{M}},
\end{equation}
and is the preneighbourhood induced from $\gamma$.

\begin{Thm}
 \label{inducedwnbdnbd}
 Let \Arr{p}{P}{X} be an admissible subobject of $X$.
 \begin{enumerate}[label=(\alph*),ref=(\alph*),align=left]
  \item \label{inducedwnbd}
	If \opair{X}{\gamma} be an internal weak neighbourhood space of
	\Bb{A} then so also is \opair{P}{\gamma_p}.

  \item \label{inducednbd}
	If \opair{X}{\gamma} be an internal neighbourhood space of
	\Bb{A} and
	\Arr{\finv{p}{}}{\Sub{X}{\mathsf{M}}}{\Sub{P}{\mathsf{M}}}
	preserve joins the \opair{P}{\gamma_p} is also an internal
	neighbourhood space of \Bb{A}.
 \end{enumerate}
\end{Thm}

\begin{proof}
 \begin{enumerate}[itemsep=1ex,label=(\alph*),align=left]
  \item Given any $m \in \Sub{P}{\mathsf{M}}$, $u \in \gamma_p(m)
	\Leftrightarrow (\exists w \in \gamma(p\circ m))\bigl(p \wedge w
	\leq p\circ u\bigr)$. Since $\gamma$ is a weak neighbourhood,
	from Theorem \ref{pnbdisnbd<>interpolative} (page
	\pageref{pnbdisnbd<>interpolative}), there exists a $v \in
	\gamma(p\circ m)$ such that $w \in \gamma(v)$. Hence
	$\finv{p}{v} \in \gamma_p(m)$ and $\finv{p}{w} \in
	\gamma_p(\finv{p}{v})$. Since $\finv{p}{w} \leq u$, $u \in
	\gamma_p(\finv{p}{v})$, showing $\gamma_p$ to be interpolative.\\

	Hence $\gamma_p$ is a weak neighbourhood on $P$. 

  \item From \ref{inducedwnbd} $\gamma_p$ is a weak neighbourhood on
	$P$. It remains to show that $\gamma_p$ preserves meets. Since
	for any family \seq{m}{i}{I} of admissible subobjects of $P$,
	$\gamma_p(\bigvee_{i \in I}m_i) \subseteq \bigcap_{i \in
	I}\gamma_p(m_i)$, it is enough to show the other
	inequality.\\

	If $u \in \bigcap_{i \in I}\gamma_p(m_i)$ then for each $i \in
	I$, there exists a $w_i \in \gamma(p\circ m_i)$ such that $p
	\wedge w_i \leq p\circ u$.\\

	Let $w = \bigvee_{i \in I}w_i$. Since \finv{p}{} preserve joins,
	using Theorem \ref{ppj full or
	partial}\ref{monoppj<>meetdistoverjoins} (page
	\pageref{monoppj<>meetdistoverjoins}) $p \wedge w = \bigvee_{i
	\in I} (p \wedge w_i) \leq p\circ u$. Hence:
	\begin{equation*}
	 w \in \bigcap_{i \in I}\gamma(p\circ m_i) = \bigcap_{i \in
	 I}\gamma(\img{p}{m_i}) = \gamma\bigl(\bigvee_{i \in
	 I}\img{p}{m_i}\bigr) = \gamma\biggl(\img{p}{\bigl(\bigvee_{i
	 \in I}m_i\bigr)}\biggr) = \gamma\biggl(p\circ\bigl(\bigvee_{i
	 \in I}m_i\bigr) \biggr),
	\end{equation*}
	implies $u \in \gamma_p\bigl(\bigvee_{i \in I}m_i\bigr)$.\\

	Hence $\gamma_p$ is a neighbourhood on $P$.
 \end{enumerate}
\end{proof}

%
%
%\begin{enumerate}[label=\underline{Remark \arabic*},ref=\underline{Remark \arabic*},align=left,resume=rem,itemsep=1.2ex]
% \item Given any preneighbourhood morphism \Arr{f}{(X, \gamma)}{(Y,
%       \phi)} and any admissible subobject \Arr{p}{P}{Y} of $Y$ it is
%       easy to see that \Arr{f_p}{(\finv{f}{P},
%       \gamma_{\finv{f}{p}})}{(P, \phi_p)} is a preneighbourhood
%       morphism and $(\finv{f}{P}, \gamma_{\finv{f}p}) = (X, \gamma)
%       \times_{(Y, \phi)} (P, \phi_p)$.
%\end{enumerate}

As expected, a regular epimorphism of a category \Bb{X} would be {\em
hereditary} if its restriction to every subobject of the codomain in
\Bb{X} is also a regular epimorphism of \Bb{X}. Explicit conditions for
$\Bb{X} = \pre{\Bb{A}}, \NHD{\Bb{A}}$ are obtained.

\subsubsection{Hereditary Regular Epimorphisms of \pre{\Bb{A}}}

\begin{Df}
 \label{hereditaryregepi}
 A regular epimorphism \Arr{f}{(X, \gamma)}{(Y, \phi)} of \pre{\Bb{A}}
 is said to be {\em hereditary} if for every admissible subobject
 \Arr{p}{P}{Y} of $Y$ the restriction $f_p$ of $f$ to $p$ in the pullback $\xymatrix{
 {\finv{f}{P}} \ar[d]_{\finv{f}{p}} \ar[r]^{f_p} & {P} \ar[d]^p \\ {X}
 \ar[r]_f & {Y} }$ is a regular epimorphism \Arr{f_p}{(\finv{f}{P},
 \gamma_{\finv{f}{p}})}{(P, \phi_p)} of \pre{\Bb{A}}. 
\end{Df}

%\begin{enumerate}[label=\underline{Remark \arabic*},
% ref=\underline{Remark\,\arabic*},align=left,resume=rem,itemsep=1.2ex] 
% \item Given any morphism \Arr{f}{X}{Y} of \Bb{A}, any admissible subobject
%       \Arr{t}{T}{Y} of $Y$ and \Arr{u}{U}{T} of $T$ the diagram:
%       \[
%       \xymatrix@!=12ex{
%       {\finv{f_t}{U}} \ar[r]^{\finv{f_t}u} \ar[d]_{(f_t)_u} & {\finv{f}{T}}
%       \ar[r]^{\finv{f}{t}} \ar[d]|{f_t} & {X} \ar[d]^f \\ 
%       {U} \ar[r]_u & {T} \ar[r]_t & {Y}
%       }
%       \]
%       suggest:
%       \begin{align}
%	\label{inv-rest-trans-eq}
%	\finv{f}(t\circ u) = \finv{f}{t}\circ \finv{f_t}{u} &
%	\quad\text{and}\quad & f_{(t\circ u)} = (f_t)_u.
%       \end{align}
%\end{enumerate}

\begin{Thm}
 \label{hereditaryregepidesc}
 A pre-neighbourhood morphism \Arr{f}{(X, \gamma)}{(Y,
 \phi)} is a hereditary regular epimorphism if and only if  for
 each $t \in \Sub{Y}{\mathsf{M}}$ \Arr{f_t}{\finv{f}T}{T} is a regular
 epimorphism of \Bb{A} and for any $u, v \in \Sub{T}{\mathsf{M}}$:
 \begin{equation}
  \label{hereditaryregepidesc-eq}
  (\exists p \in \gamma(\finv{f}{(t\circ u)}))\bigl(\img{f}{(\finv{f}t
  \wedge p)} \leq t\circ v\bigr) \Rightarrow \\
  (\exists q \in \phi(t\circ u))\bigl(t \wedge q \leq t\circ v\bigr).
  \end{equation}
\end{Thm}

\begin{proof}
 Since $\finv{f}{t\circ u} = \finv{f}{t}\circ\finv{f_t}{u}$, $p \in
 \gamma(\finv{f}{(t\circ u)}) 
 = \gamma((\finv{f}t)\circ \finv{f_t}u)$, $\img{f}{(\finv{f}t \wedge p)}
 \leq t\circ v \Leftrightarrow \finv{f}t \wedge p \leq \finv{f}{(t \circ
 v)} = (\finv{f}t)\circ(\finv{f_t}{v})$. Hence the hypothesis of
 \eqref{hereditaryregepidesc-eq} in view of \eqref{subobjectprenbd-eq}
 is equivalent to $\finv{f_t}v \in \gamma_{\finv{f}t}(\finv{f_t}u)$.

 In view of \eqref{subobjectprenbd-eq} the consequent of
 \eqref{hereditaryregepidesc-eq} is equivalent to $v \in \phi_t(u)$.

 Since \Arr{f_t}{(\finv{f}T, \gamma_{\finv{f}t})}{(T,
 \phi_t)} is a pre-neighbourhood morphism,
 \eqref{hereditaryregepidesc-eq} is equivalent to stating:
 \[
  v \in \phi_t(u) \Leftrightarrow \finv{f_t}v \in
 \gamma_{\finv{f}t}{(\finv{f_t}u)}. 
 \]

 Further, if \Arr{f}{X}{Y} is hereditarily in $\mathsf{E}$ then
 $\finv{f}v \in \gamma(\finv{f}u) \Rightarrow \finv{f}u \leq \finv{f}v
 \Leftrightarrow \img{f}{\finv{f}u} = u \leq v$.

 The equivalence now follows from the description of regular
 epimorphisms of \pre{\Bb{A}} in Theorem \ref{regepipreobj-desc} (page
 \pageref{regepipreobj-desc}). 
\end{proof}

\subsubsection{Conditions ensuring hereditary regular epimorphisms of \pre{\Bb{A}}}

The regular epimorphisms of \pre{\Set} are hereditary (see
\cite{BentleyHerrlichLowen1991}). The purpose of this and the next
subsection is to explain this phenomena. 

Consider the diagram in Figure \ref{conditions-eq1} (page
\pageref{conditions-eq1}) for the admissible subobjects $p \in
\Sub{X}{\mathsf{M}}$ and $t \in \Sub{Y}{\mathsf{M}}$. With respect to
Figure \ref{conditions-eq1}: 
\begin{itemize}
 \item the {\em blue} arrows indicate morphisms from $\mathsf{M}$ while
       the {\em orange} arrows indicate morphisms from $\mathsf{E}$,
 \item the front, right and left hand vertical squares are pullback
       squares, 
 \item the bottom horizontal square is the
       \fact{\mathsf{E}}{\mathsf{M}} of \comp{f}{p},
 \item hence the vertical squares on the right and left are completely in
       $\mathsf{M}$, \hfill{and}
 \item the diagonal on the top horizontal square is the
       \fact{\mathsf{E}}{\mathsf{M}} of
       \comp{f_t}{\bigl(\ffinv{f}{t}{p}\bigr)}. 
\end{itemize}

Since: $t\circ f_t\circ \finv{\bigl(\finv{f}t\bigr)}p = f\circ
\finv{f}t\circ \finv{\bigl(\finv{f}t\bigr)}p = f\circ p\circ
\bigl(\finv{f}t\bigr)_p = \img{f}p\circ
\rest{f}{P}\circ\bigl(\finv{f}t\bigr)_p$, it follows from the right hand
vertical pullback square the existence of a unique morphism
\Arr{w}{\finv{(\finv{f}t)}P}{\finv{t}{\img{f}{P}}} making the top
horizontal and hind vertical squares to commute.

Since the vertical diagonal with vertices
$\bigl(\ffinv{f}{t}{P}\bigr)$-$T$-$Y$-$P$ is a composite of the front
and left 
vertical pullback squares, it is a pullback square; since this is also a
composite of the hind and right hand vertical squares, and the right
hand vertical square is a pullback, it follows that the hind vertical
square is also a pullback square.

Further from the commutative square $\xymatrix@!=18ex{ {\ffinv{f}{t}{P}}
\ar[r]^{\rest{f_t}{\ffinv{f}{t}{P}}}
\ar[d]_w & {\img{f_t}{\ffinv{f}{t}{P}}} 
\ar[d]^{\img{f_t}{\ffinv{f}{t}{p}}}
\ar@{-->}[dl]|{!\, r} \\
{\finv{t}{\img{f}{P}}}  
\ar[r]_{\finv{t}{\img{f}{p}}} & {T} }$, since
$\rest{f_t}{\ffinv{f}{t}{P}} \in \mathsf{E}$ and $\finv{t}{\img{f}p} \in
\mathsf{M}$ there exists a unique morphism 
\Arr{r}{\img{f_t}{\ffinv{f}{t}{P}}}{\finv{t}{\img{f}{P}}} making the
whole diagram to commute.  

Hence, the top left triangle on the top horizontal square yields a
\fact{\mathsf{E}}{\mathsf{M}} of $w$, entailing:
\begin{align}
 \label{imagerestffinv-eq}
  \img{f}{\bigl(p \wedge \finv{f}t\bigr)} = t\circ\img{f_t}{\ffinv{f}tp} & &
  \text{ and } & & \rest{f}{(P \cap \finv{f}{T})} =
 \rest{f_t}{\ffinv{f}{t}{P}}, \\
 \intertext{and from the existence of $r$:}
 \label{img-invwedgeleqwedge-img}
 \img{f_t}{\ffinv{f}tp} \leq \finv{t}{\img{f}p} & \Leftrightarrow &
 t\circ\img{f_t}{\ffinv{f}tp} \leq \img{f}p & \Leftrightarrow & 
 \img{f}{\bigl(p \wedge \finv{f}t\bigr)} \leq t \wedge \img{f}p.
\end{align}

\begin{center}
 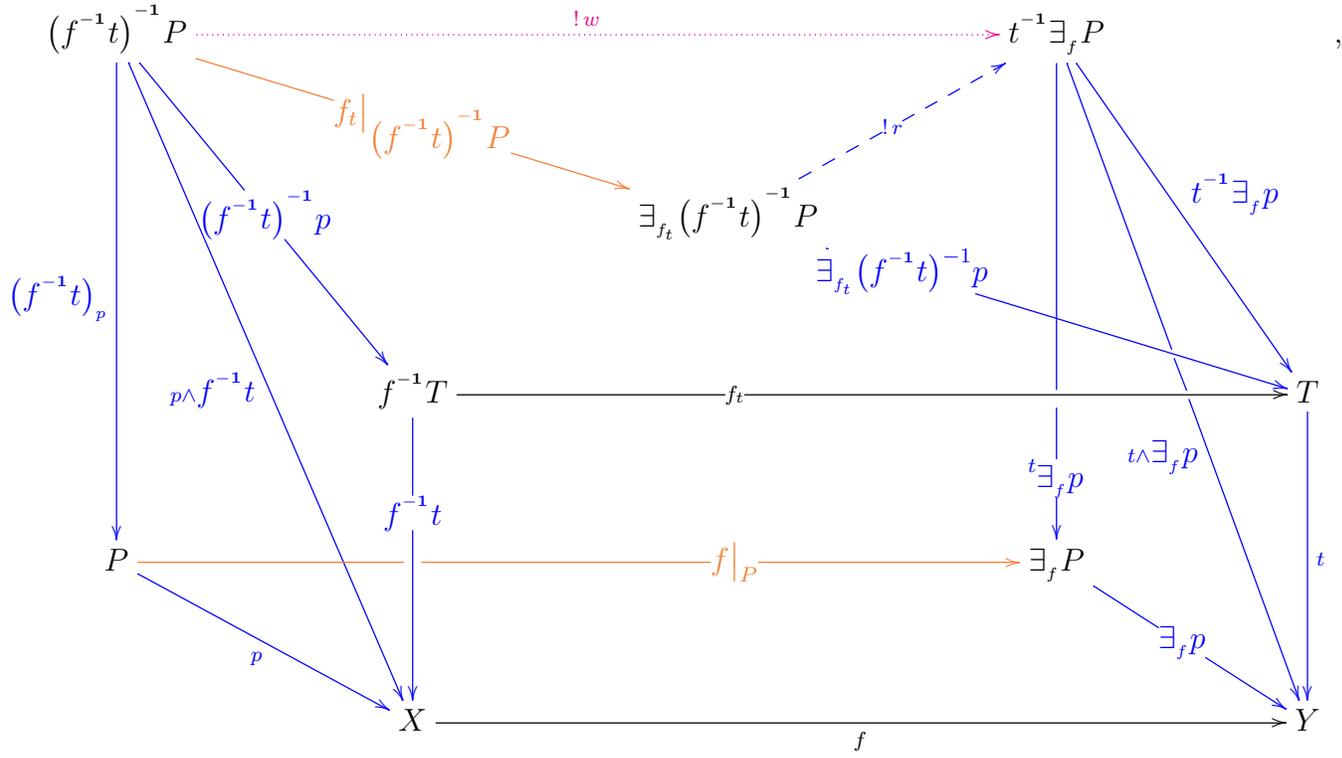
\begin{figure}[p]
  \begin{equation*}
    \xymatrixrowsep{0.6in}
    \xymatrixcolsep{0.9in}
    \xymatrix{
    {\finv{\bigl(\finv{f}t\bigr)}{P}}
    \ar@[blue][ddd]_{\color{blue} \bigl(\finv{f}t\bigr)_p}
    \ar@[blue][ddr]|{\color{blue} \finv{\bigl(\finv{f}t\bigr)}{p}}
    \ar@[Magenta]@{.>}[rrr]^{\color{Magenta}!\, w} 
    \ar@[blue][ddddr]_{\color{blue} p \wedge
    \finv{f}t} 
    \ar@[Orange][drr]|{\color{Orange}
    \rest{f_t}{\finv{\bigl(\finv{f}t\bigr)}{P}}} & & & 
    {\finv{t}{\img{f}{P}}} 
    \ar@[blue]'[dd][ddd]|{\color{blue} t_{\img{f}p}}
    \ar@[blue][ddr]^{\color{blue} \finv{t}{\img{f}p}}
    \ar@[blue][ddddr]_(0.6){\color{blue} t \wedge
    \img{f}p}|!{[dl];[ddr]}{\color{White}a} |!{[ddll];[ddr]}{\color{White}a}
    & \\
   & & {\img{f_t}{\finv{\bigl(\finv{f}t\bigr)}{P}}}
    \ar@[blue][drr]|(0.3){\color{blue}\img{f_t}{\bigl(\finv{f}t\bigr)^{-1}{p}}}
    \ar@[blue]@{-->}[ur]|{\color{blue}!\, r} & & \\ 
   & {\finv{f}{T}} \ar[rrr]|(0.36){f_t}
    \ar@[blue][dd]|(0.36){\color{blue} \finv{f}t} & & & {T}
    \ar@[blue][dd]^{\color{blue} t} \\ 
   {P} \ar@[Orange]'[r][rrr]|{\color{Orange} \rest{f}{P}}
    \ar@[blue][dr]_{\color{blue} p} & & & {\img{f}P}
    \ar@[blue][dr]|{\color{blue} \img{f}p} & \\
   & {X} \ar[rrr]_f & & & {Y}
    },
  \end{equation*}
  \caption{Frobenius morphisms}
  \label{conditions-eq1}
 \end{figure}
\end{center}

\begin{Df}
 \label{frobeniuspair-df}
 An adjunction \adjtsit{f}{g}{X}{Y} between partially ordered sets is
 said to be a {\em Frobenius pair} if:
 \begin{equation}
  \label{Frobenius-eq}
   f(g(y) \wedge x) = y \wedge f(x), \quad\text{for all }x\in X, y
   \in Y.
 \end{equation}

 If for a given morphism \Arr{f}{P}{Q} of \Bb{A} the adjunction
 \adjt{\img{f}{}}{\finv{f}{}} is a Frobenius pair then $f$ is
 a {\em Frobenius morphism}.
\end{Df}

In case of the category \Set of sets and functions every function is a
Frobenius morphism. The discussion preceding the definition above
produces equivalent formulations for Frobenius morphisms in
categories with a proper factorisation system.

\begin{Thm}
 \label{equivalenceofFrobenius}
 \tfae for any morphism \Arr{f}{X}{Y} of \Bb{A}:
 \begin{enumerate}[label=(\alph*),ref=(\alph*),align=left]
  \item \label{Frobenius} $f$ is a Frobenius morphism.
  \item \label{roundsquares} For every admissible subobject
	\Arr{t}{T}{Y} of $Y$ the diagram:
	\[
	\xymatrix{
	{\Sub{X}{\mathsf{M}}} \ar[r]^{\img{f}{}}
	\ar[d]_{\finv{\bigl(\finv{f}t\bigr)}{}} & 
	{\Sub{Y}{\mathsf{M}}} \ar[d]^{\finv{t}{}} \\ 
	{\Sub{\finv{f}T}{\mathsf{M}}} \ar[r]_{\img{f_t}{}} & 
	{\Sub{T}{\mathsf{M}}} }	 
	\]
	of order preserving maps commute.   
  \item \label{uniqueisinE} For every admissible subobject \Arr{t}{T}{Y}
	of $Y$ and \Arr{p}{P}{X} of $X$ the unique morphism
	\Arr{w}{\finv{(\finv{f}t)}P}{\finv{t}{\img{f}P}} in Figure 
	\ref{conditions-eq1} (see page \pageref{conditions-eq1}) is in
	$\mathsf{E}$. 
 \end{enumerate}
\end{Thm}

\begin{enumerate}[label=\underline{Remark \arabic*},
 ref=\underline{Remark\,\arabic*},align=left,resume=rem,itemsep=1.2ex] 
 \item The \fact{\mathsf{E}}{\mathsf{M}} system is said to satisfy {\em
       Beck-Chevalley} condition if for every pullback diagram
       $\xymatrix{ {X \times_Z Y} \ar[r]^(0.6){f_g} \ar[d]_{g_f} & {Y}
       \ar[d]^g \\ {X} \ar[r]_f & {Z} }$ the diagram:
       \[
       \xymatrixcolsep{1.2in}
       \xymatrix{
       {\Sub{X}{\mathsf{M}}} \ar[r]^{\img{f}{}} \ar[d]_{\finv{g_f}{}} &
       {\Sub{Z}{\mathsf{M}}} \ar[d]^{\finv{g}{}} \\
       {\Sub{X \times_Z Y}{\mathsf{M}}} \ar[r]_(0.6){\img{f_g}{}} &
       {\Sub{Y}{\mathsf{M}}} 
       }
       \]
       of order preserving maps commute.\\

       The diagram in \ref{roundsquares} of the Theorem is a special
       case for $g$ an admissible subobject. Hence $f$ is a Frobenius
       morphism if and only if  the \fact{\mathsf{E}}{\mathsf{M}}
       system satisfies Beck-Chevalley condition for admissible
       subobjects of codomain of $f$. 
\end{enumerate}

\subsubsection{Five Conditions for Heredity of Regular Epimorphisms of \pre{\Bb{A}}}

The case for regular epimorphisms for \pre{\Set} to be hereditary is a
consequence of every function being a Frobenius morphism. The
heredity of regular epimorphisms in \pre{\Bb{A}} holds in fact with
weaker conditions as the following theorem establishes.

\begin{Thm}
 \label{condensurherregepi}
 Given the statements:
 \begin{enumerate}[label=(\alph*),ref=(\alph*),align=left]
  \item \label{stableE} The set $\mathsf{E}$ is stable under pullbacks.
  \item \label{hereditaryE} The set $\mathsf{E}$ is hereditary.
  \item \label{Frobenuis} Every morphism of \Bb{A} is Frobenuis.
  \item \label{FrobenuisE} Every morphism in $\mathsf{E}$ is Frobenius.
  \item \label{Frobenuisregepi} Every regular epimorphism is Frobenius.
  \item \label{regepihereditary} Every regular epimorphism of
	\pre{\Bb{A}} is hereditary.
 \end{enumerate}

 The following implications hold good:
 \[
 \xymatrix{
 {\ref{stableE}} \ar@2{->}[r] & {\ref{hereditaryE}} \ar@2{<->}[d] & & \\
 & {\ref{Frobenuis}} \ar@2{->}[r] & {\ref{FrobenuisE}} \ar@2{->}[r] &
 {\ref{Frobenuisregepi}} \ar@2{->}[r] & {\ref{regepihereditary}}
 }.
 \]
\end{Thm}

\begin{proof}
 Follows immediately from the diagram in Figure \ref{conditions-eq1}
 (page \pageref{conditions-eq1}) and Theorem
 \ref{equivalenceofFrobenius} (page
 \pageref{equivalenceofFrobenius}). The equivalence of \ref{hereditaryE}
 and \ref{Frobenuis} is obvious.  
\end{proof}

\begin{enumerate}[label=\underline{Remark \arabic*},
 ref=\underline{Remark\,\arabic*},align=left,resume=rem,itemsep=1.2ex] 
 \item The obvious equivalence of
 \ref{hereditaryE} and \ref{Frobenuis} in Theorem
       \ref{condensurherregepi} was also observed
       in \cite[Proposition 1.3]{ClementinoGiuliTholen1996}.

 \item It is known from \cite{JanelTholen1994} that the condition
       \ref{stableE} in Theorem \ref{condensurherregepi} is equivalent
       to the \fact{\mathsf{E}}{\mathsf{M}} system satisfying the
       Beck-Chevalley condition.
\end{enumerate}

In the case of \Set, since $\mathsf{E} = \mathrm{Epi}$ is pullback
stable, every regular epimorphism of the category \pre{\Set} of
pretopological spaces is hereditary.

\subsubsection{Hereditary Regular Epimorphisms of \NHD{\Bb{A}}}

Since the morphisms of internal neighbourhood spaces have the
preimage preserve join property, in view of Theorem
\ref{inducedwnbdnbd}\ref{inducednbd} (page \pageref{inducednbd}) it is
best to restrict to the case when every morphism of \Bb{A} has preimage
preserve join property. Hence, from Corollary \ref{ppj>subobjectframes}
(page \pageref{ppj>subobjectframes}), every lattice of admissible
subobjects is a frame.

\begin{Df}
 \label{hereditaryregepinbd}
 A regular epimorphism \Arr{f}{(X, \gamma)}{(Y, \phi)} of \NHD{\Bb{A}}
 is said to be {\em hereditary} if for every admissible subobject
 \Arr{p}{P}{Y} of $Y$ the restriction $f_p$ of $f$ to $p$ in the
 pullback $\xymatrix{ {\finv{f}{P}} \ar[d]_{\finv{f}{p}} \ar[r]^{f_p} &
 {P} \ar[d]^p \\ {X} \ar[r]_f & {Y} }$ is a regular epimorphism
 \Arr{f_p}{(\finv{f}{P}, \gamma_{\finv{f}{p}})}{(P, \phi_p)} of
 \pre{\Bb{A}}.  
\end{Df}

\begin{enumerate}[label=\underline{Remark \arabic*},
 ref=\underline{Remark\,\arabic*},align=left,resume=rem,itemsep=1.2ex]
 \item If \Arr{f}{(X, \gamma)}{(Y, \phi)} is a hereditary regular
       epimorphism of \NHD{\Bb{A}}, then each restriction
       \Arr{f_p}{(\finv{f}{P}, \gamma_{\finv{f}{p}})}{(Y, \phi_p)} is a
       regular epimorphism of \NHD{\Bb{A}}. Hence each $f_p \in
       \RegEpi{\Bb{A}} \subseteq \mathsf{E}$ and each \finv{f_p}{}
       preserves arbitrary joins.\\

       Furthermore, for each $p \in \Sub{Y}{\mathsf{M}}$,
       $\img{f}{\finv{f}{p}} = p$.
    
 \item \label{leq<>finvleq}
       For each $p, q \in \Sub{Y}{\mathsf{M}}$:
       \[
	\finv{f}{p} \leq \finv{f}{q} \Leftrightarrow
       \img{f}{\finv{f}{p}} \leq q \Rightarrow p \leq q.
       \]

       Hence:
       \[
	p \leq q \Leftrightarrow \finv{f}{p} \leq \finv{f}{q}.
       \]

       Since for any $u \leq p \in \Sub{Y}{\mathsf{M}}$,
       $\finv{f}{p}\circ \finv{f_p}{u} = \finv{f}{(p\circ u)}$, the same
       holds for each $f_p$.
       
 \item \label{restrnbd-alt}
       Consequently, for each $p \in \Sub{Y}{\mathsf{M}}$:
       \[
	u \in \phi_p(m) \Leftrightarrow \finv{f_p}{u} \in
       \gamma_{\finv{f}{p}}(\finv{f_p}{m}). 
       \]
\end{enumerate}

\begin{Thm}
 \label{hereditaryregepiinnbd}
 Assume every morphism of \Bb{A} has the preimage preserve join
 property and \Arr{f}{(X, \gamma)}{(Y, \phi)} is a morphism of
 \NHD{\Bb{A}}. 
 
 \begin{enumerate}[label=(\alph*),ref=(\alph*),align=left,itemsep=1.2ex]
  \item \label{herregepinbd>fsurregepipre}
	If $f$ is a hereditary regular epimorphism of \NHD{\Bb{A}} then
	it is regular epimorphism of \pre{\Bb{A}} such that for each $p
	\in \Sub{Y}{\mathsf{M}}$, $f_p \in \RegEpi{\Bb{A}}$. 

  \item \label{Frobregepipre>herregepinbd}
	If $f$ is a Frobenius morphism, $f_p \in \RegEpi{\Bb{A}}$ for
	each $p \in \Sub{Y}{\mathsf{M}}$ and a regular epimorphism of
	\pre{\Bb{A}} then it is a hereditary regular epimorphism of
	\NHD{\Bb{A}}. 

  \item \label{regepipre<>pseudoopen}
	If $f$ is a regular epimorphism of \Bb{A} with the property that
	for each $p \in \Sub{Y}{\mathsf{M}}$, $f_p \in \mathsf{E}$, then
	$f$ is a regular epimorphism of \pre{\Bb{A}}, if
	and only if, for every $y \in \Sub{Y}{\mathsf{M}}$: 
	\begin{equation}
	 \label{pseudopen-eq}
	 u \in \gamma(\finv{f}{y}) \Rightarrow \img{f}{u} \in \phi(y).
	\end{equation}
 \end{enumerate}
\end{Thm}

\begin{proof}
 \begin{enumerate}[label=(\alph*),ref=(\alph*),align=left,itemsep=1.2ex]
  \item Follows from \ref{leq<>finvleq}, Theorem \ref{regepinbdobj-desc}
	(page \pageref{regepinbdobj-desc}) \& Theorem
	\ref{regepipreobj-desc} (page \pageref{regepipreobj-desc}).

  \item It is required to show for any $p \in \Sub{Y}{\mathsf{M}}$ the
	morphism \Arr{f_p}{(X, \gamma_{\finv{f}{p}})}{(Y, \phi_p)} is a
	regular epimorphism of \NHD{\Bb{A}}.

	Given $u, m, p \in \Sub{Y}{\mathsf{M}}$, $u, m \leq p$, in
	light of \ref{restrnbd-alt}, it is enough to show:
	\[
 	 \finv{f_p}{u} \in \gamma_{\finv{f}{p}}{\finv{f_p}{m}}
	\Rightarrow u \in \phi_p(m).
	\]

	Clearly, from equation \eqref{subobjectprenbd-eq} (page
	\pageref{subobjectprenbd-eq}):
	\begin{multline*}
	 \finv{f_p}{u} \in \gamma_{\finv{f}{p}}(\finv{f_p}{m})
	 \\
	 \Leftrightarrow \biggl(\exists t \in
	 \gamma\bigl((\finv{f}{p})\circ (\finv{f_p}{m})\bigr)\biggr)
	 \bigl(\finv{f}{p} \wedge t \leq
	 (\finv{f}{p})\circ(\finv{f_p}{u})\bigr) \\ 
	 \Leftrightarrow \biggl(\exists t \in
	 \gamma\bigl(\finv{f}{(p\circ m)}\bigr)\biggr)
	 \bigl(\finv{f}{p} \wedge t \leq \finv{f}{(p \circ u)}\bigr) \\
	 \Leftrightarrow \biggl(\exists t \in
	 \gamma\bigl(\finv{f}{(p\circ m)}\bigr)\biggr)
	 \bigl(\img{f}{(\finv{f}{p} \wedge t)} \leq p\circ u\bigr) \\
	 \Leftrightarrow \biggl(\exists t \in
	 \gamma\bigl(\finv{f}{(p\circ m)}\bigr)\biggr)
	 \bigl(p \wedge \img{f}{t} \leq p\circ u\bigr) \quad\text{
	 (since $f$ is Frobenius)} \\
	 \Rightarrow (\exists t \in
	 \Sub{X}{\mathsf{M}})\biggl(\finv{f}{\img{f}{t}} \in
	 \gamma(\finv{f}{(p\circ m)}) \text{ and } p \wedge \img{f}{t}
	 \leq p \circ u\biggr) \quad (\text{ since
	 }\adjt{\img{f}{}}{\finv{f}{}}) \\
	 \Rightarrow (\exists t \in
	 \Sub{X}{\mathsf{M}})\biggl(\img{f}{t} \in \phi(p\circ m) \text{
	 and } p \wedge \img{f}{t} \leq p \circ u\biggr) \\
	 \Rightarrow u \in \phi_p(m),
	\end{multline*}
	completing the proof of \ref{Frobregepipre>herregepinbd}.

  \item If \Arr{f}{(X, \gamma)}{(Y, \phi)} is a morphism of
	\NHD{\Bb{A}} then 
	\adjt{\img{f}{}}{\adjt{\finv{f}{}}{\fInv{f}{}}}, where for any
	$t \in \Sub{X}{\mathsf{M}}$:
	\[
	 \img{f}{t} = \bigwedge\bigl\{p \in \Sub{Y}{\mathsf{M}}: t \leq
	\finv{f}{p}\bigr\} \quad\text{ and }\quad
	 \fInv{f}{t} = \bigvee\bigl\{q \in \Sub{Y}{\mathsf{M}}:
	\finv{f}{q} \leq t\bigr\}.
	\]

	Hence for any $t \in \Sub{X}{\mathsf{M}}$:
	\begin{equation}
	 \label{cond-forall<=exists}
	 \fInv{f}{t} \leq \img{f}{t} \Leftrightarrow (\forall p, q \in
	\Sub{Y}{\mathsf{M}})\bigl(\finv{f}{q} \leq t \leq \finv{f}{p}
	\Rightarrow q \leq p\bigr).	  
	\end{equation}

	Since for each $p \in \Sub{Y}{\mathsf{M}}$, $f_p \in
	\mathsf{E}$, and \ref{leq<>finvleq} shows the statement on
	the right hand side of \eqref{cond-forall<=exists} holds,
	and hence for all $t \in \Sub{X}{\mathsf{M}}$, $\fInv{f}{t} \leq
	\img{f}{t}$.

	Assume now for each $y \in \Sub{Y}{\mathsf{M}}$:
	\[
	 u \in \gamma(\finv{f}{y}) \Rightarrow \img{f}{u} \in \phi(y). 
	\]

	If $\finv{f}{p} \in \gamma{\finv{f}{y}}$ then $p =
	\img{f}{\finv{f}{p}} \in \phi(y)$, showing $f$ to be a regular
	epimorphism of \pre{\Bb{A}} using \ref{restrnbd-alt} (page
	\pageref{restrnbd-alt}) and Theorem \ref{regepipreobj-desc}
	(page \pageref{regepipreobj-desc}).

	Conversely if $f$ be a regular epimorphism of \pre{\Bb{A}} then:
	\[
	u \in \gamma(\finv{f}{y}) \Rightarrow \finv{f}{y} \leq u
	\Leftrightarrow y \leq \fInv{f}{u} \Rightarrow y \leq
	\fInv{f}{u} \leq \img{f}{u},
	\]
	and $u \leq \finv{f}{\img{f}{u}} \Rightarrow
	\finv{f}{\img{f}{u}} \in \gamma(\finv{f}{y})$ implies
	$\img{f}{u} \in \phi(y)$.
 \end{enumerate}
\end{proof}

\begin{enumerate}[label=\underline{Remark \arabic*},
 ref=\underline{Remark\,\arabic*},align=left,resume=rem,itemsep=1.2ex]
 \item \label{pseudoopendef}
       In the case of \Set, neighbourhood morphisms \Arr{f}{(X,
       \gamma)}{(Y, \phi)} satisfying \eqref{pseudopen-eq} are called
       {\em pseudo open}. \\ 

       We could call a morphism of \NHD{\Bb{A}} {\em pseudo open} if it
       satisfies \eqref{pseudopen-eq}. \\

       Hence: a morphism of \NHD{\Bb{A}} is a regular epimorphism of
       \pre{\Bb{A}} if and only if  it is a pseudo open regular
       epimorphism of \Bb{A} with each of its restrictions in
       $\mathsf{E}$. 
\end{enumerate}
\section{Concluding Remarks}
\label{Examples}

\subsection{The category \Set}\label{sets}
The concepts studied in this paper are well known for the category \Set
of sets and functions.

\Set comes equipped with its usual \fact{Epi}{Mono} system. The lattice
\Sub{X}{Mono} of admissible subobjects is a complete atomic Boolean
algebra. Hence, $\Int{\Top}{\Set} = \Top = \NHD{\Set}$ is the category
of topological spaces. The category \pre{\Set} is isomorphic to the
category $\mathbf{preTop}$ of pretopological spaces. The category
$\mathbf{preTop}$ is investigated in \cite{BentleyHerrlichLowen1991} \&
\cite{HerrlichLowenSchwarz1991}.

Since epimorphisms in \Set are pullback stable the regular epimorphisms
of pretopological spaces are hereditary (see \cite[Theorem 26, page
14]{BentleyHerrlichLowen1991} and compare Theorem
\ref{condensurherregepi}, page \pageref{condensurherregepi}). However,
it is also known from 
\cite[Example 4, page 5]{BentleyHerrlichLowen1991}, that the regular epimorphisms 
of pretopological spaces are not in general pullback stable.

\subsubsection{A Weak Neighbourhood which is not a Neighbourhood}
Neighbourhoods in \Set can be obtained by just specifying the filters
for each point, since the subobject lattices are atomic. This is not
true of weak neighbourhoods or preneighbourhoods.

Given a set $X$ and a topology $\Theta$ on $X$ let \pow{\Theta}{c} be
the set of closed subsets of the topological space
\opair{X}{\Theta}. Define: 
\[
 \mu(M) = \bigl\{V \subseteq X: (\exists C \in \Theta^c)(M \subseteq C
 \subseteq V)\bigr\}. 
\]

Clearly, $\mu$ defines a preneighbourhood on $X$ such that $\mo_\mu =
\Theta^c$, and:
\begin{align*}
 \intr{\mu}{M} = \bigcup\bigl\{C \in \Theta^c: C \subseteq M\bigr\}, \\
 \mu(M) = \bigcup\bigl\{\atleast{C}: C \in \Theta^c\text{ and }M
 \subseteq C\bigr\}. 
\end{align*}

Hence $\mu$ is a weak neighbourhood, and a neighbourhood if and only
if \pow{\Theta}{c} is closed under arbitrary joins. Incidentally, under
the same condition \intr{\mu}{} becomes a Kurastowski interior.

\subsection{The category \Top}\label{topologies}

The category \Top of topological spaces comes equipped with its usual
\fact{Epi}{ExtMon} system. The lattice \Sub{X}{ExtMon} is precisely the
set of all subsets of $X$ equipped with the subspace topology and hence
again is a complete atomic Boolean algebra.

A preneighbourhood \Arr{\CAL{F}}{\opp{\Sub{X}{ExtMon}}}{\Fil{X}} on a
topological space $X$ is given on specifying for 
each $T \subseteq X$ a filter $\CAL{F}_T$ of subspaces of $X$ such that $S
\in \CAL{F}_T \Rightarrow T \subseteq S$. Thus, for instance, taking all
open sets (or, closed sets) containing $T$ provides instances of two
preneighbourhood structures on $X$.  

Since neighbourhoods are meet
preserving and the subobject lattices are atomic, it is enough to
specify the neighbourhoods of each $x \in X$. Thus, neighbourhoods on
$X$ correspond to specifying a second topology on $X$. Consequently,
\NHD{\Top} is isomorphic to the category \BiTop of bitopological spaces
and functions which are continuous with respect to both the topologies
on $X$.

\subsection{The category \Loc}\label{locales}
The category \Loc of locales comes equipped with its usual
\fact{Epi}{RegMon} system. The lattice \Sub{X}{RegMon} is a distributive
complete lattice in which finite joins distribute over arbitrary meets,
i.e., is a coframe. For any localic map \Arr{f}{X}{Y} the preimage is
\Arr{f_{-1}}{\Sub{Y}{RegMon}}{\Sub{X}{RegMon}}, defined as the largest
sublocale of $X$ which is contained inside the subset \finv{f}{S} ($S
\in \Sub{Y}{RegMon}$) of $X$ (see \cite[Chapter III.4.2, page
29]{PicadoPultr2012}). It is known for every localic map $f$ its
preimage $f_{-1}$ preserve, apart from arbitrary meets, finite joins
(see \cite[Theorem III.9.2, page 41]{PicadoPultr2012}). 

Recently in \cite{DubeIghedo2016} \& \cite{DubeIghedo2016a}
neighbourhoods have been effectively used. The neighbourhood used is
\Arr{\mathit{o}_X}{\opp{\Sub{X}{RegMon}}}{\Fil{X}}:
\begin{equation}
 \label{naturaltopologyonlocales-eq}
 \mathit{o}_X(S) = \bigl\{T \in \Sub{X}{RegMon}: (\exists a \in X)(S \subseteq
 \slopen{a} \subseteq T)\bigr\}. 
\end{equation}

Since \Arr{\mathfrak{o}}{X}{\Sub{X}{RegMon}} preserves finite meets
and arbitrary joins, for any $b \in X$ and any family \seq{a}{i}{I} of
elements of $X$:
\begin{multline*}
 \mathfrak{o}(b) \cap \bigvee_{i \in I}\mathfrak{o}(a_i) =
 \mathfrak{o}(b) \cap \mathfrak{o}\bigl(\bigvee_{i \in I}a_i\bigr) =
 \mathfrak{o}\bigl(b \wedge \bigvee_{i \in I}a_i\bigr) \\
 = \mathfrak{o}\bigl(\bigvee_{i \in I}(b \wedge a_i)\bigr) = \bigvee_{i
 \in I}\slopen{b \wedge a_i} = \bigvee_{i \in I}\bigl(\slopen{b} \cap
 \slopen{a_i}\bigr), 
\end{multline*}
shows $\mathbf{OpenSub}(X)$ the set of open sublocales of $X$ is a
frame. Further, since $\slopen{a} \leq \slopen{b} \Leftrightarrow a \leq
b$, \Arr{\mathfrak{o}}{X}{\mathbf{OpenSub}(X)} is an isomorphism
of frames.  

Since $\mo_{\mathit{o}_{X}} = \mathbf{OpenSub}(X)$, it follows
\opair{X}{\mathit{o}} is actually an internal topological space of
\Loc. 

Furthermore, for any frame homomorphism \Arr{f}{X}{Y}, if $S \in
\Sub{Y}{RegMon}$ is a sublocale of $Y$ and $T \in \mathit{o}_Y(S)$, then
there exists a $b \in 
Y$ such that $S \subseteq \slopen{b} \subseteq T$. Hence:
\[
 f_{-1}S \subseteq f_{-1}\slopen{b} = \slopen{f^*(b)} \subseteq f_{-1}T,
\]
where \Arr{f^*}{Y}{X} is the left adjoint of $f$, which is a frame
homomorphism. This implies $f_{-1}T \in \mathit{o}_X(f_{-1}S)$,
yielding: 
\begin{Thm}
 \label{Thembanbdrtinv}
 The functor \Arr{\mathit{O}}{\Loc}{\pNHD{\Loc}} defined by
 $\mathit{O}(X) = (X, \mathit{o}_X)$ is a right inverse to the forgetful
 functor \Arr{U}{\pNHD{\Loc}}{\Loc}. 
\end{Thm}

\subsection{Acknowledgments}\label{acknowledgements}
I am indebted to:
\begin{enumerate}
 \item T. Dube for supporting my research through NRF Funds
       for Research Chair here at Unisa.
 \item M. Korostenski-Davies for painstakingly going through the draft
       version of this document and suggesting several editorial
       changes.
 \item Z. Janelidze for his stimulating ideas during our talks on
       several occasions.
 \item A. Razafindrakato, D. Holgate and their students for sharing
       their experiences.
\end{enumerate}

\afterpage{\clearpage}

%\input{Conclusions}

%backmatter

%\printindex

%%% insert the right number of ../'s from the current directory
%\bibliographystyle{ams}
%\bibliography{../../../tex/essentials/bib-algebra,../../../tex/essentials/bib-category,../../../tex/essentials/bib-frames,../../../tex/essentials/bib-topology}

\printbibliography

\end{document}